%% file: erosion_methods.tex
\documentclass[preprint, 10pt]{elsarticle}

\input{preambles.tex}
\begin{document}

\title{A boundary-integral framework to simulate \\ viscous erosion of a porous medium}


\author[Bryan]{Bryan D.~Quaife}
\author[Nick]{M.~Nicholas J.~Moore}
\address[Nick]{Department of Mathematics and Geophysical Fluid Dynamics Institute, Florida State University, Tallahassee, FL, 32306.}
\address[Bryan]{Department of Scientific Computing and Geophysical Fluid Dynamics Institute, Florida State University, Tallahassee, FL, 32306.}

\begin{abstract} 
We develop numerical methods to simulate the fluid-mechanical erosion of many bodies in two-dimensional Stokes flow. The broad aim is to simulate the erosion of a porous medium (e.g.~groundwater flow) with grain-scale resolution. Our fluid solver is based on a second-kind boundary integral formulation of the Stokes equations that is discretized with a spectrally-accurate Nystr\"om method and solved with fast-multipole-accelerated GMRES.  The fluid solver provides the surface shear stress which is used to advance solid boundaries. We regularize interface evolution via curvature penalization using the {\thL} formulation, which affords numerically stable treatment of stiff terms and therefore permits large time steps.  The overall accuracy of our method is spectral in space and second-order in time. The method is computationally efficient, with the fluid solver requiring $\mathcal{O}(N)$ operations per GMRES iteration, a mesh-independent number of GMRES iterations, and a one-time $\mathcal{O}(N^2)$ computation to compute the shear stress. We benchmark single-body results against analytical predictions for the limiting morphology and vanishing rate. Multibody simulations reveal the spontaneous formation of channels between bodies of close initial proximity. The channelization is associated with a dramatic reduction in the resistance of the porous medium, much more than would be expected from the reduction in grain size alone.
\end{abstract}


\maketitle

\section{Introduction\label{s:intro}}
Flow-induced erosion deteriorates and reshapes solid material over a range of scales found in nature, from massive land formations~\cite{han1969, perron2009formation, abrams2009growth, petroff2011geometry, jerolmack2012internal, Rothman2012, coh-dev-sey-yi-szy-rot2015, perkins2015amplification}, to centimeter-scale features and patterns~\cite{daerr2003erosion, berhanu2012shape, nienhuis2014wavelength, domokos2014river}. Though less visible, these same forces are working at the very smallest scales, slowly deteriorating the individual constituents of porous media (e.g.~soil, sand, or clay) or biological structures like plaque and biofilms~\cite{pic-van-hei2000, sha2002, gro-gij-van-fer-hat-van-yua-wen2007}. Motivated by such examples, this paper develops numerical methods to simulate fluid-mechanical erosion of multiple solid bodies in Stokes flow---the most relevant regime for groundwater \cite{bear2013dynamics, cao2010coupled} and the aforementioned bio-fluid applications.

Recent work in the high-Reynolds-number regime has highlighted the importance of a shape-flow feedback that occurs during erosion~\cite{ris-moo-chi-she-zha2012, moo-ris-chi-zha-she2013, hewett2017evolution, moore2017riemann, lopez2018cfd} and the related processes of dissolution and melting~\cite{Huang2015, kondratiuk2015steadily, dallaston2015channelization, rycroft2016asymmetric, cohen2016erosion, hewett2017pear, claudin2017dissolution}.  In these cases, the fluid alters the morphology of immersed structures, which in turn modifies the surrounding flow field, and so on. Such examples of `sculpting with flow'~\cite{ristroph2018sculpting} have been examined in the context of a single body~\cite{ris-moo-chi-she-zha2012, moore2017riemann} and, more recently, multiple bodies whose shape evolution is intertwined~\cite{hewett2017evolution}. From this viewpoint, the example of an eroding porous medium presents new opportunities and challenges. First, since the medium is composed of many individual constituents, multiple bodies must be considered. Second, the dense packing present in typical porous-media applications requires close-range interactions to be resolved. Finally, long-range effects inherent to the Stokes limit couples the shape evolution of {\em all} bodies to a greater degree than any other regime.

\begin{figure}
\begin{center}
\includegraphics[width = 0.80 \textwidth]{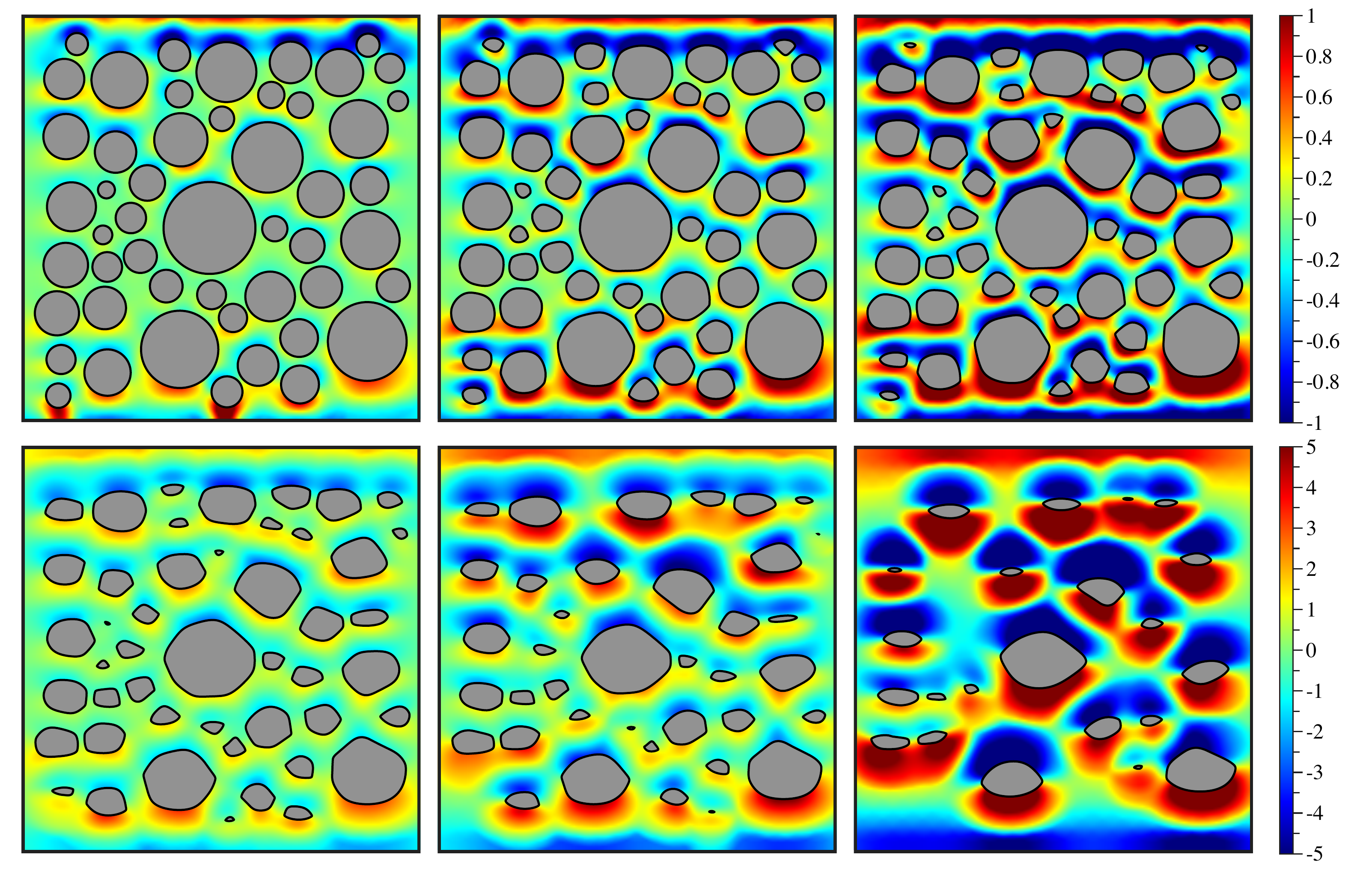}
\caption{\label{fig:50bodies} Simulation of 50 bodies eroding in Stokes flow under the action of fluid-mechanical shear. The flow is horizontal (left to right) and the 6 snapshots are evenly spaced in time. Color represents vorticity, which provides a convenient way to visualize local shear rates. Erosion not only diminishes the size of the bodies but also alters their shapes considerably. A few well-defined channels develop in the space between bodies. Several of the bodies vanish in finite time and the simulation continues without interruption.
}
\end{center}
\end{figure}

Ultimately, the methods developed here will lead to high-fidelity simulations of multibody Stokes erosion as illustrated in Fig.~\ref{fig:50bodies} for the case of 50 bodies (of randomly assigned initial sizes and positions). Here, the erosion of solid material arises from the mechanical stress induced by the surrounding flow. In particular, we take the local erosion rate to be proportional to the magnitude of the shear stress---a law that has been verified by laboratory experiments~\cite{ris-moo-chi-she-zha2012}. As seen in Fig.~\ref{fig:50bodies}, erosion not only diminishes the size of the bodies, but also alters their shapes considerably. In particular, several bodies develop corner-like features and nearly flat faces, with straight channels intervening between them. In the figure, color represents the vorticity of the surrounding flow, which reduces to shear on solid boundaries and thus provides a convenient way to visualize local erosion rates. Several bodies vanish in finite time and the simulation continues without interruption.

\paragraph{Contributions}
To perform the above simulation, two fast and accurate methods are
employed---one to solve the fluid equations and the other to evolve the
eroding bodies.  The fluid equations are reformulated in terms of a boundary integral equation (BIE)~\cite{poz1992}.  Since the boundary conditions are Dirichlet, we employ a double-layer potential so that the resulting BIE is second-kind.  This second-kind BIE is discretized with a spectrally accurate trapezoid rule~\cite{tre-wei2014} which results in a well-conditioned dense linear system for an unknown density function. The linear system is solved using fast-multipole-accelerated GMRES~\cite{saa-sch1986}.  Once the density function is determined, we compute the shear stress with a spectrally accurate odd-even quadrature rule~\cite{sid-isr1988}.

At the next stage, we must use the shear stress to advance solid boundaries. Tracking fluid-structure interfaces with long time horizons presents numerical challenges in that the underlying mesh can become distorted and tangled. Furthermore, the process of erosion naturally leads to regions of high curvature and even corners.  We address these tracking issues in two ways: first, we introduce a tangential velocity that keeps the grid equispaced with respect to arclength but that does not affect dynamics; second, we suppress high-frequency oscillations in the interface evolution via curvature penalization. The interface is parametrized in the so-called {\thL} coordinate system so that the curvature penalization appears as a diffusive term \cite{hou-low-she1994}. This term can be easily be handled in a numerically stable fashion with an exponential integrator.  Once the simulation is complete, we perform a post-processing step in which we compute bulk properties, such as vorticity and pressure, using an upsampled trapezoid rule.  

We supplement our primarily computational study with a few analytical results. In particular, our simulation shows that a single eroding body tends to a slender, for-aft symmetric form characterized by nearly uniform shear stress along its surface. This body exhibits sharp corners at its front and rear. The principle of long-time uniform stress combined with some local analysis leads to an analytical prediction for the opening angles formed at the front and rear. Additionally, slender-body theory gives rise to scaling laws for the reduction rate of the body's area and drag. These analytical results, in addition to their own intrinsic value, allow us to conclusively validate the numerical simulation.

\paragraph{Limitations}
We make two fundamental design choices in the simulations described here: first, we restrict attention to two spatial dimensions, and second the bodies are considered immobile (i.e.~fixed to their initial locations). By restricting to two dimensions, we can simulate many bodies (up to $50$ in this manuscript) with high resolution, enabling accurate characterization of porous-medium erosion. While BIEs are available in three dimensions, the extra computational cost limits the number of bodies that could be practically simulated and makes porous-medium simulations less feasible. Meanwhile, the second assumption of immobile bodies affords a simple framework to explore the effects of erosion in isolation. While other mechanisms, such as sedimentation or transport, undoubtedly play a role in real applications, these effects would introduce additional complexity, both conceptually and numerically. We therefore defer their study to future work.

\paragraph{Related Work}
In a closely related study, Mitchell \& Spagnolie used a traction integral equation formulation to simulate the viscous erosion of a single body in three dimensions~\cite{mit-spa2016}. The flexibility of this formulation enabled simulations of a {\em mobile} body eroding in an {\em arbitrary} background flow (e.g.~uniform flow, shear flow, etc.).  An advantage of the traction formulation is that it directly provides the shear stress needed to evolve surfaces.  In the present work, we aim to simulate the erosion of a porous medium and, further, choose to fix the bodies to their initial locations. Thus, rather than a {\em single, mobile} body as considered in~\cite{mit-spa2016}, we aim to simulate the erosion of {\em many, immobile} bodies. Although the traction formulation has several attractive features, we instead employ a standard double-layer formulation since it allows us to easily calculate bulk quantities (e.g.~vorticity, pressure, and velocity) that will assist our understanding of the erosion process.

For the porous network simulations considered here, it is essential to
resolve nonlocal effects between the different bodies.  BIEs naturally
resolve these interactions, and they have been used extensively to
simulate viscous fluids in complex geometries, including porous media
flow~\cite{dea-qua-bir-jua2018, bar-mar-vee-zha2018}.  Common challenges
of a BIE formulation include developing efficient preconditioners for
the discretized system and numerical difficulties associated with nearly
touching bodies. In this work, we consider geometries requiring
$\mathcal{O}(100)$ GMRES iterations, so preconditioners are currently
unnecessary. In future work, we will address mobile bodies for which
near contacts may necessitate preconditioners~\cite{qua-bir2015a, qua-cou-dar2018, cou-pou-dar2017}. Similarly, near-singular schemes will not be necessary due to the fact that the spacing between receding bodies of fixed position always grows in time. Again, the future consideration of mobile bodies may necessitate near-singular schemes~\cite{qua-bir2014a, klo-bar-gre-one2013, bar-wu-vee2015, hel-oja2008a, bea-lai2001}.

\paragraph{Outline of the Paper}
In Section~\ref{s:formulation}, we describe the governing equations for both the fluid and the erosion model in a \thL~framework.  In Section~\ref{s:method}, the numerical methods in both space and time are described and Section~\ref{s:qoi} describes how to compute the vorticity, pressure, and drag.  Section~\ref{s:SingleResults} gives several numerical examples for a single body, and Section~\ref{s:MultiResults} presents numerical examples for multiply-connected domains.  Finally, concluding remarks are made in Section~\ref{s:conclusions}.

\section{Formulation}
\label{s:formulation}
\begin{figure}[htpb]
  \centering
  \input{schematic.tikz}
  \caption{\label{fig:schematic} A schematic for the governing equations.  A no-slip boundary condition is imposed on each body $\gamma_\ell$ whose unit normal points outward relative to the geometry.  On the outer geometry $\Gamma$, a Hagen-Poiseuille flow is imposed.  The bodies erode at a rate that is equal to the magnitude of the shear stress applied by the fluid.  The bodies are constrained to the middle third of the channel that is located between the dashed lines.}
\end{figure}
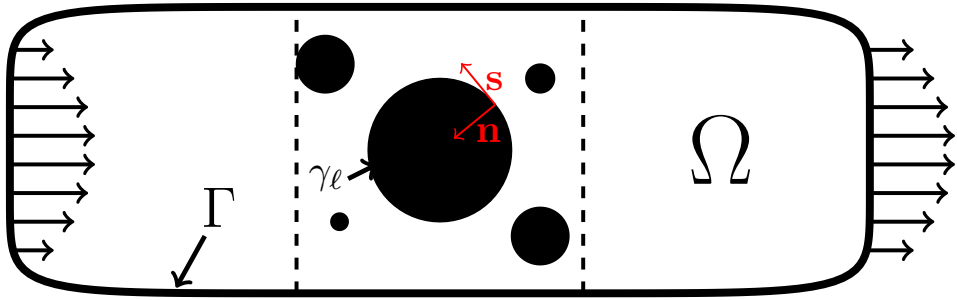

We start by defining the main variables used to model erosion.  We consider flows that are confined by a solid wall $\Gamma$ that encloses $M$ eroding bodies as diagrammed in Fig.~\ref{fig:schematic}.  The bodies are denoted as $\gamma_\ell$, $\ell=1,\ldots,M$, and we write $\gamma = \gamma_1 \cup \cdots \cup \gamma_M$.  Therefore, the boundary of the fluid domain is $\bd\Omega = \Gamma \cup \gamma$. On each body, $\nn$ and $\ss$ denote the unit normal and unit tangent vectors respectively.  Throughout the paper, we adopt the convention that $\ss$ points in the counterclockwise direction and $\nn$ points out of the fluid domain (i.e.~into the bodies). Neglecting inertial forces, the dynamics of the fluid is fully determined by the position of the bodies $\xx_\ell(s,t) \in \gamma_\ell$, where $s$ is the arclength and $t$ is time.  Given $\xx_\ell(s,t)$, $\ell=1,\ldots,M$, derived variables include the fluid velocity $\uu$, the pressure $p$, and the shear stress $\tau$.   On the bounding wall, we prescribe a velocity $\UU(\xx,t)$.  Then, the governing equations are
\begin{equation}
\label{eqn:erosionModel}
\begin{split}
  \mu \Delta \uu = \grad p, &\hspace{20pt} \xx \in \Omega, \gap &&\mbox{conservation
of momentum}\\
\grad \cdot \uu = 0, &\hspace{20pt} \xx \in \Omega, \gap &&\mbox{\em conservation of mass} \\
\uu = 0, &\hspace{20pt} \xx \in \gamma, \gap &&\mbox{\em no slip on the
bodies} \\
\uu = \UU, &\hspace{20pt} \xx \in \Gamma, \gap &&\mbox{\em outer wall
velocity} \\
\Vn = \abs{\tau}, &\hspace{20pt} \xx \in \gamma,
&&\mbox{\em erosion model},
\end{split}
\end{equation}
where $\Vn$ is the normal velocity of the interface and $\mu$ is the fluid viscosity.  A scaling constant often appears in the erosion model in~\eqref{eqn:erosionModel}, but in the Stokes limit, the constant only sets the time scale, so we have taken it to be one.  For the sake of numerical stability, we will modify the erosion model in Section~\ref{sec:thetaL}. 

The outer wall is formed by rounding off the corners of $[-3,3] \times [-1,1]$, and we impose the Hagen-Poiseuille flow
\begin{align}
  \UU(\xx) = \umax (1-y^2,0)\, , \quad \xx \in \Gamma\, ,
\end{align}
where $U$ sets the maximum of the imposed Poiseuille flow. Often, we will simply set $U=1$, but in some instances we change $U$ dynamically to enforce that the pressure drop across the computational cell is constant as bodies erode within it.  To minimize the boundary effects of the inflow and outflow, we only place bodies in the center third, $[-1,1] \times [-1,1]$, of the fluid domain.  

We describe how the incompressible Stokes equations are solved in $\Omega$ using a boundary integral equation in Section~\ref{sec:bies}.  In Section~\ref{sec:shearStressLP}, we discuss how the shear stress is computed.  Finally, in Section~\ref{sec:thetaL}, we describe how corners in the bodies are avoided and perform a change of coordinates to treat a stiff term.

\subsection{Boundary integral equation formulation} 
\label{sec:bies}

In this section, we describe an equivalent BIE formulation of the incompressible Stokes equations with Dirichlet boundary conditions
\begin{equation}
  \begin{split}
  \mu\Delta \uu &= \grad p\, , \qquad && \xx \in \Omega, \\
  \grad \cdot \uu &= 0,   && \xx \in \Omega, \\
  \uu &= \ff,  && \xx \in \bd\Omega.
  \end{split}
  \label{eqn:stokes}
\end{equation}
The boundary condition $\ff$ is the prescribed velocity $\UU$ on
$\Gamma$, and it is zero on the eroding bodies $\gamma$.  In two dimensions, a BIE for a stream function formulation is possible~\cite{gre-kro-may1996}, but this does not extend to three dimensions.  Therefore, we adopt a primitive variable formulation where the velocity is written in terms of the standard double-layer potential~\cite{lad1963,poz1992}
\begin{align}
  \DD[\eeta](\xx) = \frac{1}{\pi} \int_{\bd\Omega} 
    \frac{\rr \cdot \nn}{\rho^2} \frac{\rr \otimes \rr}{\rho^2} 
    \eeta(\yy) \, ds_\yy\, , \quad \xx \in \Omega \, ,
    \label{eqn:velocityDLP}
\end{align}
where $\rr = \xx - \yy$, $\rho = \|\rr\|$, and $\eeta$ is an unknown density function.  As the viscosity $\mu$ also sets the time scale, we take $\mu=1$.

Since the double-layer potential is not able to represent force- and torque-free velocity fields,  we complete the double-layer potential with Stokeslets and rotlets~\cite{pow-mir1987, pow1993}
\begin{align}
  S[\llambda_\ell,\cc_\ell] = \frac{1}{4\pi} \left(-\log \rho \II + 
    \frac{\rr \otimes \rr}{\rho^2} \right)\llambda_\ell
  \quad \text{and} \quad
  R[\xi_\ell,\cc_\ell] = \frac{\rr^\perp}{\rho^2}\xi_\ell\, , 
  \quad \ell=1,\ldots,M \, ,
  \label{eqn:stokeslet_rotlet}
\end{align}
where  $\rr = \xx - \cc_\ell$, $\cc_\ell$ is a point inside the $\ell^{th}$ grain, $\rr^\perp = (r_2,-r_1)$, and $\rho = \|\rr\|$.  Then, the solution of~\eqref{eqn:stokes} is
\begin{align}
  \label{eqn:completed_DLP}
  \uu(\xx) = \DD[\eeta](\xx) + 
    \sum_{\ell=1}^{M} S[\llambda_\ell,\cc_\ell](\xx) +
    \sum_{\ell=1}^{M} R[\xi_\ell,\cc_\ell](\xx)\, , \quad \xx \in \Omega \, ,
\end{align}
where the density function, Stokeslets, and rotlets satisfy
\begin{subequations}
\label{eqn:completed_BIE}
\begin{align}
  \ff(\xx) &= -\frac{1}{2}\eeta(\xx) + \DD[\eeta](\xx) + 
      \sum_{\ell=1}^{M} S[\llambda_\ell,\cc_\ell](\xx) +
      \sum_{\ell=1}^{M} R[\xi_\ell,\cc_\ell](\xx) + \NN_0[\eeta](\xx)\, ,
      \quad \xx \in \bd\Omega \, , 
      \label{eqn:DLP}\\
  \llambda_\ell &= \frac{1}{2\pi}\int_{\gamma_\ell} \eeta(\yy) \,
  ds_\yy\, ,
  \quad \ell=1,\ldots,M \, ,
  \label{eqn:stokeslet} \\
  \xi_\ell &= \frac{1}{2\pi}\int_{\gamma_\ell} \yy^\perp \cdot \eeta(\yy)
  \, ds_\yy\, , \quad \ell=1,\ldots,M \, .
  \label{eqn:rotlet}
\end{align}
\end{subequations}
Here,
\begin{align}
  \NN_0[\eeta](\xx) = \int_{\Gamma} 
    (\nn(\xx) \otimes \nn(\yy))\eeta(\yy) \, ds_\yy, \quad \xx \in
    \Gamma\, ,
\end{align}
removes the rank one null space resulting from the flux-free requirement of $\ff$.

This BIE formulation requires that $\bd\Omega$ is smooth.  Unfortunately, smoothness is not guaranteed by erosion since corners develop at locations where the shear stress vanishes.  One option is to modify our BIE to handle corners~\cite{rac-ser2017, ser-rok2016, hel2011, gil-hao-mar2014, bre2012}, but these methods require the corner's location, a non-uniform spaced mesh, or are ill-conditioned.  Instead, we slightly modify the erosion model in Section~\ref{sec:thetaL} so that we can control the smoothness of the bodies.  In this manner, the BIE as described above can be applied.

\subsection{Computing the shear stress}
\label{sec:shearStressLP}
Once the density function, Stokeslets, and rotlets are computed, any desired quantity, including the shear stress, can be computed.  The shear stress is
\begin{align}
  \tau = -\left(\nabla \uu + \nabla \uu^T \right)\nn \cdot \ss \, .
\end{align}
The reader may notice an unusual sign convention above, but this is simply a consequence of our conventions for $\nn$ and $\ss$ (see Fig.~\ref{fig:schematic}).  The velocity field of the completed double-layer potential~\eqref{eqn:completed_DLP} includes velocities from three terms: the double-layer potential, the Stokeslets, and the rotlets.  We compute the deformation tensor $\ssigma = \frac{1}{2}(\nabla\uu + \nabla\uu^T)$ of these terms individually.

The deformation tensor of the double-layer potential~\eqref{eqn:velocityDLP} at an interior point $\xx \in \Omega$ is~\cite{qua-bir2014a}
\begin{equation}
  \label{eqn:shearStressDLP}
  \begin{aligned}
  \ssigma^\DD(\xx) = \frac{1}{2\pi}\int_{\bd\Omega} &\left(
    2\frac{\rr \cdot \nn}{\rho^2} \frac{\rr \cdot \eeta}{\rho^2} \II + 
    \frac{\rr \cdot \eeta}{\rho^4} (\nn \otimes \rr + \rr \otimes \nn) 
    \right. \\
    &\left.
    +\frac{\rr \cdot \nn}{\rho^4} (\eeta \otimes \rr + \rr \otimes \eeta) - 
    8\frac{(\rr \cdot \nn)(\rr \cdot \eeta)}{\rho^6}(\rr \otimes \rr)
  \right) \eeta(\yy) \, ds_\yy \, .
  \end{aligned}
\end{equation}
We require the deformation tensor on $\bd\Omega$, and this requires taking the limit of equation~\eqref{eqn:shearStressDLP} when $\xx$ tends to $\xx_0 \in \bd\Omega$.  The limiting value is
\begin{align}
  \lim_{\substack{\xx \rightarrow \xx_0 \\ \xx \in \Omega}}\sigma_\DD(\xx) =
  J[\eeta](\xx_0) + \ssigma^\DD(\xx_0)\, , \quad \xx_0 \in \bd\Omega \, ,
\end{align} 
where
\begin{align}
  J[\eeta](\xx_0) = \frac{1}{2}\left(\pderiv{\eeta}{\ss} 
    \cdot \ss\right) \left[ 
  \begin{array}{cc}
    s_x^2 - s_y^2 & 2s_x s_y \\ 2s_x s_y & s_y^2 - s_x^2
  \end{array}\right] \, .
  \label{eqn:shearStressJump}
\end{align}
Next, the deformation tensor of the Stokeslets and rotlets~\eqref{eqn:stokeslet_rotlet} are
\begin{equation}
  \label{eqn:shearStressSR}
  \begin{aligned}
  \ssigma^S(\xx_0) &= \sum_{\ell=1}^{M}
    \frac{\rr \cdot \llambda_\ell}{4\pi\rho^2} \left(
    \II - \frac{2}{\rho^2} \rr \otimes \rr \right),  \\
  \ssigma^R(\xx_0) &= -\sum_{\ell=1}^M
    \frac{\xi_\ell}{\rho^4} \left(\rr \otimes \rr^\perp + 
    \rr^\perp \otimes \rr \right) \, .
  \end{aligned}
\end{equation}
Combining these deformation tensors, the shear stress at $\xx_0 \in \bd\Omega$ is
\begin{align}
  \tau = -2 \left(J[\eeta](\xx_0) + \ssigma^\DD(\xx_0) + 
    \ssigma^S(\xx_0) + \ssigma^R(\xx_0)\right) \nn \cdot \ss \, .
  \label{eqn:shearStressLP}
\end{align}

\subsection{Boundary evolution in the {\thL} framework} 
\label{sec:thetaL}

To evolve the boundaries of eroding bodies, we use the {\thL} framework which offers certain advantages in the ability to numerically stabilize moving-interface problems~\cite{hou-low-she1994}. Here, $\theta$ is the local tangent angle and $L$ is the total perimeter of a given boundary. For a boundary parametrized by arclength $(x(s),y(s))$, the local tangent angle is defined through the relation $\left( \pdi{x}{s}, \pdi{y}{s} \right) = \left(\cos \theta, \sin \theta \right)$.  

For numerical stability, we modify the erosion law~\eqref{eqn:erosionModel} to include a term that penalizes curvature $\kappa$. The idea is to prevent regions of extremely high curvature from developing, as these regions would compromise the accuracy and stability of our fluid solver. We thus replace~\eqref{eqn:erosionModel} with
\begin{align}
  \label{eqn:Vn}
  \Vn = \atau + \eps \mean{\atau} \left( \frac{L }{2\pi} \kappa - 1
  \right) \, .
\end{align}
The parameter $\eps \ll 1$ controls the strength of the curvature-penalization. This term causes regions of high curvature to recede more rapidly and thus suppresses any high-frequency oscillations that may exist. Above, $\mean{\cdot}$ indicates a boundary mean. Since a closed body has mean curvature $\mean{\kappa} = 2\pi / L$, the penalty term has mean zero and thus preserves area. The only source of material loss is therefore the shear-stress itself. In addition, we have scaled the penalty term on the mean absolute shear, $\mean{\atau}$, so that the ratio of smoothing-to-erosion always remains fixed, even as bodies vanish.  

In the following, we parameterize each boundary by normalized arclength, $\alpha = 2 \pi s / L$. In terms of $\alpha$, the curvature is given by
\begin{align}
\kappa = \frac{2 \pi}{L} \pderiv{\theta}{\alpha} \, ,
\end{align}
and substitution into~\eqref{eqn:Vn} gives
\begin{align}
\label{eqn:Vn2}
  \Vn = \atau +  \eps \mean{\atau}   \left(\pderiv{\theta}{\alpha} - 1 \right) \, .
\end{align}

The boundary of each body will be represented by a grid of collocation points. If we were to simply evolve the interfaces under~\eqref{eqn:Vn2}, this grid would become distorted, which, if unchecked, could lead to severe tangling and numerical instability. We therefore introduce an artificial tangential velocity $\Vs$ to keep the grid points equally spaced in arclength~\cite{hou-low-she1994}. The appropriate tangential velocity satisfies 
\begin{align}
  \pderiv{\Vs}{\alpha} = \thalpha \Vn - \mean{\thalpha \Vn} \, .
\end{align} 
The evolution of the interface is then given by the combination of the normal and tangential velocities
\begin{align}
  \dot{\xx}(t) = \Vn \nn + \Vs \ss \, .
\end{align}
Since $\Vs$ is tangential, it does not modify the physics of erosion, only the parameterization of the boundaries.

When recast in the variables $\theta$ and $L$, the interface evolution becomes
\begin{align}
  \tderiv{L}{t} &= - 2 \pi \mean{\thalpha \Vn}, \\
  \pderiv{\theta}{t} &= \frac{2 \pi}{L} \left( \pderiv{\Vn}{\alpha} + \thalpha \Vs \right) \, .
\end{align}
Inserting~\eqref{eqn:Vn2} for the normal velocity gives
\begin{align}
\label{eqn:Levo1}
\tderiv{L}{t} &= - 2 \pi \mean{\thalpha \Vn}, \\
\label{eqn:thetaevo1}
  \pderiv{\theta}{t} &= \frac{2\pi}{L} \left(
  \eps \mean{\atau} \ppd{\theta}{\alpha} + \pderiv{\atau}{\alpha} +
  \thalpha \Vs \right)\,.
\end{align}
These are the main evolution equations for $\theta(\alpha,t)$ and $L(t)$.  In addition, it is necessary to track a reference point for each body, so that the {\thL} variables can be converted back to Cartesian coordinates. We track the coordinates of the surface mean $\mean{\xx}$, which moves according to
\begin{align}
\label{eq:xsurf}
\frac{d}{dt} \mean{\xx} = \mean{ V_s \ss + V_n \nn } \, .
\end{align}
We note that tracking the center of mass is also a valid choice, but the surface mean gives a somewhat simpler formula. Equations~\eqref{eqn:Levo1}--\eqref{eq:xsurf} comprise the main system of equations to be solved numerically in Section~\ref{sec:timeStepping}. 

\section{Numerical methods}
\label{s:method}
To perform numerical simulations, we exploit the difference between the erosion and fluid time scales.  Since erosion occurs much slower, we freeze the geometry at each time step, solve the incompressible Stokes equations, compute the shear stress, and then update the geometry according to the erosion model.  Therefore, there are three main numerical methods that need to be developed.  In Section~\ref{sec:BIE}, we describe how the BIE~\eqref{eqn:completed_BIE} is solved.  Once the density function is found, Section~\ref{sec:shearStress} describes how the shear stress~\eqref{eqn:shearStressLP} is computed.  Then, Section~\ref{sec:timeStepping} describes a time stepping method for advancing the interface according to equations~\eqref{eqn:Levo1}--\eqref{eq:xsurf}. 

To create stable simulations with long time horizons, we use numerical methods that achieve spectral accuracy in space, second-order accuracy in time, and maintain an equispaced discretization of the bodies.  The boundaries are represented in Fourier space, giving us a set of discretization points $\xx_i$ of each body, and all derivatives and integrals (the arclength, curvature, and mean values) are computed with spectral accuracy using the Fourier representation.

\subsection{Solving the BIE}
\label{sec:BIE}

To numerically compute the completed double-layer potential~\eqref{eqn:completed_BIE}, we discretize each of the $M$ bodies with  $N_\iin$ points, and the outer boundary with $N_\out$ points.  Then, $N = MN_\iin + N_\out$ is the total number of discretization points of $\bd\Omega$.  The double-layer potential in~\eqref{eqn:DLP} is discretized with the trapezoid rule as
\begin{align}
  \DD[\eeta](\xx_i) = \sum_{j=1}^{N} K(\xx_i,\xx_j) \eeta(\xx_j) 
      \Delta s_{j}\, , \quad i=1,\ldots,N \, ,
  \label{eqn:trapDLP}
\end{align}
where $\Delta s_j$ is the arclength term of $\bd\Omega$ at
$\xx_j$, and
\begin{align}
  K(\xx,\yy) = \frac{1}{\pi} \frac{\rr \cdot \nn}{\rho^2} 
      \frac{\rr \otimes \rr}{\rho^2} \, ,
\end{align}
where $\rr = \xx - \yy$ and $\rho = \|\rr\|$, is the kernel of the double-layer potential.  The diagonal term of~\eqref{eqn:trapDLP} is replaced with the limiting value of the double-layer potential
\begin{align}
  D(\xx,\xx) = -\frac{\kappa(\xx)}{2\pi}(\ss(\xx) \otimes \ss(\xx)) \, ,
\end{align}
where $\kappa(\xx)$ is the curvature and $\ss(\xx)$ is the unit tangent vector at $\xx \in \bd\Omega$.  The equations for the Stokeslets~\eqref{eqn:stokeslet} and the rotlets~\eqref{eqn:rotlet} are also discretized with the trapezoid rule.  Applying this discretization strategy to~\eqref{eqn:completed_BIE}, the dense linear system to be solved is
\begin{subequations}
  \label{eqn:completed_BIE_discretized}
  \begin{align}
  \label{eqn:discretized_DLP}
    \ff(\xx_i) &= -\frac{1}{2}\ssigma(\xx_i) + \sum_{j=1}^{N} 
      K(\xx_i,\xx_j) \eeta(\xx_j) \Delta s_{j} + 
      \sum_{\ell=1}^M S[\llambda_\ell,\cc_\ell](\xx_i) +
      \sum_{\ell=1}^M R[\xi_\ell,\cc_\ell](\xx_i) \\
    \llambda_\ell &= \sum_{\xx_j \in \gamma_\ell} \eeta(\xx_j) 
      \Delta s_j \, , \\ 
    \xi_\ell &= \sum_{\xx_j \in \gamma_\ell}
      {(\xx_j - \cc_\ell)}^\perp \cdot \eeta(\xx_j) \Delta s_j \, ,
  \end{align}
\end{subequations}
for $i=1,\ldots,N$, $\ell=1,\ldots,M$.  When $\xx_i \in \Gamma$, the term 
\begin{align}
  \sum_{\xx_j \in \Gamma} (\nn(\xx_i) \otimes
      \nn(\xx_j))\eeta(\xx_j) \, ,
\end{align}
is added to~\eqref{eqn:discretized_DLP}.  Since we discretize integrands that are smooth and periodic, the trapezoid rule employed above achieves spectral accuracy~\cite{tre-wei2014}.  

The linear system~\eqref{eqn:completed_BIE_discretized} is solved iteratively with the generalized minimal residual method (GMRES)~\cite{saa-sch1986}.  Since the BIE is second-kind, only a mesh-independent number of GMRES iterations are required~\cite{cam-ips-kel-mey-xue1996}.  Therefore, the asymptotic cost is determined by the cost of performing a matrix-vector multiplication, and the bulk of the cost is evaluating the $N$-term sum in~\eqref{eqn:discretized_DLP}. We use the Fast Multipole Method (FMM)~\cite{gre-rok1987, gre-gre-may1992} so that the cost of solving~\eqref{eqn:completed_BIE_discretized} is $\bigO(N)$ operations.
 
\subsection{Shear stress}
\label{sec:shearStress}
We have decomposed the deformation tensor into contributions from three different terms.  The contributions from the Stokeslets and the rotlets~\eqref{eqn:shearStressSR} can be immediately evaluated since they do not involve any integration to approximate.  For the shear stress of the double-layer potential, we use a slight variant of the trapezoid rule. The kernel of the deformation tensor due to the double-layer potential~\eqref{eqn:shearStressDLP} is 
\begin{equation}
\begin{aligned}
  K_\ssigma(\xx,\yy)\eeta(\yy) &= \frac{1}{2\pi} \left(
    2\frac{\rr \cdot \nn}{\rho^2} \frac{\rr \cdot \eeta}{\rho^2} \II + 
    \frac{\rr \cdot \eeta}{\rho^4} (\nn \otimes \rr + \rr \otimes \nn)
    \right. \\ & \hspace{25pt}+ \left.
    \frac{\rr \cdot \nn}{\rho^4} (\eeta \otimes \rr + \rr \otimes \eeta) - 
    8\frac{(\rr \cdot \nn)(\rr \cdot \eeta)}{\rho^6}(\rr \otimes \rr)
  \right) \, .
\end{aligned}
\end{equation}
This kernel has a singularity that is asymptotic to $\bigO(\rho^{-2})$ when $\yy \rightarrow \xx$, and this is too strong for the unmodified trapezoid rule.  The singularity strength can be reduced since the velocity corresponding to a constant density function is constant~\cite{poz1992}.  Therefore, the deformation tensor of the double-layer potential is equivalent to
\begin{align}
  \ssigma^{\DD}(\xx) = \int_{\bd\Omega}K_\ssigma(\xx,\yy)
      (\eeta(\yy) - \eeta(\xx)) \, ds_{\yy}\, , \quad \xx \in \bd\Omega \, .
\end{align}
This expression has a singularity that is asymptotic to $\bigO(\rho^{-1})$ when $\yy \rightarrow \xx$.  Then, to achieve spectral accuracy, we apply the trapezoid rule with odd-even integration~\cite{sid-isr1988}.  Therefore, the deformation tensor of the double-layer potential at a discretization point $\xx_i$ on body $\gamma_\ell$, with $i$ being odd, is approximated as
\begin{align}
  \ssigma^{\DD}(\xx_i) = \sum_{\xx_j \notin \gamma_\ell}
    K_\ssigma(\xx_i,\xx_j) \eeta(\xx_j) \Delta s_j + 
  \sum_{\substack{\xx_j \in \gamma_\ell \\ j \: \mathtt{even}}}
    K_\ssigma(\xx_i,\xx_j) (\eeta(\xx_j) - \eeta(\xx_i)) \Delta s_j \, ,
  \label{eqn:stressOddEven}
\end{align}
and a similar approximation is used when $i$ is even.  We also need to approximate the jump in the deformation tensor~\eqref{eqn:shearStressJump}, and this is done with Fourier differentiation.  Finally, the shear stress is computed by multiplying the deformation tensor by the normal vector, and taking the inner product with the tangent vector.

We do not employ an FMM to evaluate the shear stress at $\xx_i$, $i=1,\ldots,N$, so the cost of evaluating the shear stress is $\bigO(N^2)$.  However, unlike the velocity double-layer potential that is evaluated at each GMRES iteration, the shear stress is only computed once per time step. In future work, we plan to use a kernel-independent FMM~\cite{yin-bir-zor2004} to evaluate the shear stress with $\bigO(N)$ operations.

\subsection{Time-stepping with the {\thL} method} 
\label{sec:timeStepping}

With the shear stress computed, we are now ready to numerically evolve the boundary of each body using the {\thL} framework. The main equations to be solved are~\eqref{eqn:Levo1}--\eqref{eq:xsurf}.  In Section~\ref{sec:thetaL} we introduced a curvature-penalizing term into the erosion law~\eqref{eqn:Vn}. This term suppresses high-frequency oscillations, but also introduces {\em stiffness} into evolution equation~\eqref{eqn:thetaevo1}. We will therefore use a combined exponential-integrator/Runge-Kutta method to solve system~\eqref{eqn:Levo1}--\eqref{eq:xsurf}. The exponential integrator will handle the stiff term in a numerically stable fashion, while the Runge-Kutta method will treat nonlinear terms explicitly. Our overall time-stepping method will be second-order accurate.

First, we must address another potential source of numerical instability. The erosion law~\eqref{eqn:Vn} depends on the {\em absolute value} of shear stress, $\atau$, which exhibits a corner anywhere $\tau$ changes sign, i.e.~at a stagnation point. The corner in $\atau$ introduces high frequencies into~\eqref{eqn:thetaevo1} which could result in either numerical instability or an extremely stringent step-size restriction. To avoid these difficulties, we smooth $\atau$ via a narrow Gaussian filter of width $\sigma \ll 1$
\begin{align}
\atausig = \int_{\gamma} \frac{1}{\sqrt{2\pi} \sigma}
 \exp \left(- \frac{\alpha'^2}{2 \sigma^2}\right) \abs{ \tau(\alpha - \alpha') } \, d\alpha' \, ,
\end{align}
with corresponding normal and tangential interface velocities
\begin{align}
& \Vnsig = \atausig +  \eps \mean{\atausig}
\left(\pderiv{\theta}{\alpha} - 1 \right), \\
& \pderiv{\Vssig}{\alpha} = \thalpha \Vnsig - \mean{\thalpha \Vnsig} \, .
\end{align}
This smoothing further suppresses the feedback of high frequencies between $\atau$ and the body shape.

We now consider the system of evolution equations given by~\eqref{eqn:Levo1}--\eqref{eq:xsurf}. For convenience, we introduce the following variables
\begin{align}
& \Mterm = - 2 \pi \mean{\thalpha \Vnsig} \, , \\
& \NLterm = \frac{2 \pi}{L} \left( \pderiv{\atausig}{\alpha} + \thalpha
\Vssig \right) \, , \\
& \elfun = \frac{2 \pi}{L}  \mean{\atausig} \, .
\end{align}
Then,~\eqref{eqn:Levo1} and~\eqref{eqn:thetaevo1} can be compactly expressed as
\begin{align}
\label{eqn:Levo2}
& \tderiv{L}{t} = \Mterm \, , \\
\label{eqn:thetaevo2}
& \pderiv{\theta}{t} = \eps \elfun \ppd{\theta}{\alpha} + \NLterm \,.
\end{align}
The diffusive term on the right of~\eqref{eqn:thetaevo2} results from the curvature penalization and is the stiffest term present. This term is multiplied by $\zeta$, which is spatially constant but depends on time (due to $L$ and $\mean{\atau}$). Meanwhile, $\NLterm$ represents all remaining terms in the evolution equation, including contributions from the shear-dependent erosion law and from the artificial tangential velocity, $\Vs$. Since the local shear stress depends on interactions between bodies, $\NLterm$ is both nonlinear and nonlocal.

We will evolve $\theta(\alpha,t)$ in spectral space and so introduce the Fourier series for $\theta$ and $\NLterm$
\begin{align}
\label{thetaSeries}
 \theta(\alpha,t) &= \alpha + \FourierSum \thhat_k e^{ i k \alpha} \, , \\
 \NLterm(\alpha,t)  &= \FourierSum \widehat{\NLterm}_k e^{ i k \alpha} \, .
\end{align}
The careful reader will notice that in~\eqref{thetaSeries} we have decomposed $\theta(\alpha)$ into a linear component and a periodic component. In spectral space,~\eqref{eqn:thetaevo2} produces the system of ODEs
\begin{align}
\tderiv{\thhat_k}{t} +  \eps k^2  \elfun \thhat_k =
\widehat{\NLterm}_k\, ,
\qquad k=-N_\iin/2,\ldots,N_\iin/2-1\, .
\end{align}
These ODEs can be simplified by applying the integrating factor
\begin{align}
\label{eq:mu}
\mu_k = \exp \left( \eps k^2 \int \elfun \, dt \right) \, ,
\end{align}
which leads to the coupled system
\begin{align}
\label{eqn:Levo3}
& \tderiv{L}{t} = \Mterm \, , \\
\label{eqn:thetaevo3}
& \tderiv{}{t}\left( \mu_k \thhat_k \right) = \mu_k \widehat{\NLterm}_k \, ,
\qquad k = -N_\iin/2, \ldots, N_\iin/2 -1 \, .
\end{align}
Note that the diffusive term from~\eqref{eqn:thetaevo2} does not appear explicitly in~\eqref{eqn:thetaevo3}, but rather is implied by the integrating factor~\eqref{eq:mu}. This is the idea of an exponential integrator. We solve the resulting system~\eqref{eq:mu}--\eqref{eqn:thetaevo3} with an explicit, second-order Runge-Kutta method (RK2). In particular, we use the {\em midpoint} RK2 method because it maximizes the smoothing effect of the diffusive term.

At the first stage of RK2 (the half step) we discretize~\eqref{eq:mu}--\eqref{eqn:thetaevo3} as
\begin{align}
\label{eq:mu12d}
& \mu_k^{(n+1/2)} = \mu_k^{(n)} \exp \left( \eps k^2 \frac{\Dt}{2}
\zeta^{(n)} \right) \, , \\
\label{eq:L12d}
& \frac{L^{(n+1/2)} - L^{(n)}}{\Dt/2} = \Mterm^{(n)} \, , \\
\label{eq:th12d}
& \frac{ \mu_k^{(n+1/2)} \thhat_k^{(n+1/2)} - \mu_k^{(n)} \thhat_k^{(n)}}{\Dt/2} 
= \mu_k^{(n)} \hat{\NLterm}_k^{(n)} \, ,
\end{align}
where the superscript indicates the time step. At the next stage of RK2 it will be necessary to evaluate derivatives at the half step $n+1/2$, thus requiring a call to the Stokes solver for which the body shapes must be known. We therefore insert~\eqref{eq:mu12d} into~\eqref{eq:th12d} and solve for $\thhat_k^{(n+1/2)}$, giving
\begin{align}
\label{eq:th12}
\thhat_k^{(n+1/2)} = \left( \thhat_k^{(n)} + \frac{1}{2} \Dt \, \hat{\NLterm}_k^{(n)} \right)
\exp \left( -\epsilon k^2 \frac{\Dt}{2} \elfun^{(n)}\right) \, .
\end{align}
Since $\zeta >0$, the exponential term always {\em suppresses} high-frequencies (exponentially with $k^2$). In other words, this is a smoothing term. With the geometry at step $n+1/2$ known, a call to the Stokes solver allows the computation of $\Mterm^{(n+1/2)}$, $\widehat{\NLterm}_k^{(n+1/2)}$, and $\zeta^{(n+1/2)}$.

At the second stage of RK2, we discretize~\eqref{eq:mu}--\eqref{eqn:thetaevo3} as 
\begin{align}
\label{eq:mu1d}
& \mu_k^{(n+1)} = \mu_k^{(n)} \exp \left( \eps k^2 {\Dt} \,
\zeta^{(n+1/2)} \right) \, , \\
\label{eq:L1d}
& \frac{L^{(n+1)} - L^{(n)}}{\Dt} = \Mterm^{(n+1/2)} \, , \\
\label{eq:th1d}
& \frac{ \mu_k^{(n+1)} \thhat_k^{(n+1)} - \mu_k^{(n)} \thhat_k^{(n)}}{\Dt} =
\mu_k^{(n+1/2)} \hat{\NLterm}_k^{(n+1/2)} \, .
\end{align}
We must eventually compute the geometry at step $n+1$, and so we insert~\eqref{eq:mu1d} into~\eqref{eq:th1d} and solve for $\thhat_k^{(n+1)}$, giving
\begin{align}
\thhat_k^{(n+1)} =  \thhat_k^{(n)} \exp \left( - \epsilon k^2 \Dt
\elfun^{(n+1/2)} \right) + \Dt \, \hat{\NLterm}_k^{(n+1/2)} \exp \left(
- \frac{1}{2} \eps k^2 \Dt \left( 2 \elfun^{(n+1/2)} - \elfun^{(n)}
\right) \right) \, .
\end{align}
Once again, both exponential terms act to smooth the geometry. Finally, to track the reference point, $\mean{\xx}$, we simply apply the same RK2 method to~\eqref{eq:xsurf}.

We point out that most exponential integrators use some form of quadrature to treat~\eqref{eq:mu}, but, here, we have simply applied the same RK2 method used for the rest of the system. We believe this to be a simple, systematic approach to arrive at a self-consistent discretization of desired accuracy because it avoids choices regarding which quadrature to use (midpoint, trapezoid, Simpson's, etc.) and can be easily generalized to other time-stepping methods, such as fourth-order Runge-Kutta.

\section{Post processing quantities of interest}
\label{s:qoi}
We further characterize the flow by computing the vorticity (Section~\ref{sec:vorticity}), pressure (Section~\ref{sec:pressure}), and drag (Section~\ref{sec:drag}).  Once equation~\eqref{eqn:completed_DLP} is solved, the vorticity and pressure will be computed at points both inside the fluid bulk and on the along individual bodies, and then the drag of individual bodies will be computed.  For target points in the fluid domain, the integrands of the layer potentials are periodic and smooth, and the trapezoid rule guarantees spectral accuracy~\cite{tre-wei2014}.  However, if the target point is close to $\bd \Omega$, a nearly-singular integrand must be integrated, and this requires upsampling that is described in Section~\ref{sec:NSI}.  When the target points is on the boundary, the integrands involve singularities that are evaluated by combining singularity subtraction and odd-even integration~\cite{sid-isr1988}.  

\subsection{Vorticity}
\label{sec:vorticity}
We start by computing the vorticity $\omega(\xx) = v_x - u_y$ for $\xx \in \Omega$ which we decompose into contributions from the double-layer potential, Stokeslets, and rotlets.  The vorticity of the double-layer potential~\eqref{eqn:velocityDLP} at $\xx \in \Omega$ is
\begin{align}
  \omega^{\DD}(\xx) = -\frac{1}{\pi}\int_{\bd\Omega} 
    \frac{(\rr \cdot \nn^\perp) + (\rr \cdot \nn)}{\rho^4}
    (\rr \cdot \eeta) \, ds_{\yy}\, .
  \label{eqn:vorticityDLP}
\end{align}
To compute~\eqref{eqn:vorticityDLP} numerically, we apply the trapezoid rule with an upsampling scheme described in Section~\ref{sec:NSI}.  Next, the vorticity due to the Stokeslet is
\begin{align}
  \omega^S(\xx) = -\frac{1}{\pi} \sum_{\ell=1}^{M} 
    \frac{\rr \cdot \llambda_\ell^\perp}{\rho^2}\, ,
\end{align}
and the vorticity of a rotlet is zero.  Therefore, the vorticity at $\xx
\in \Omega$ is
\begin{align}
  \omega(\xx) = \omega^\DD(\xx) + \omega^S(\xx)\, .
\end{align}

For the Stokes equations, on boundaries that have a no-slip boundary condition, the shear stress is equivalent to the vorticity.  Therefore, the vorticity is computed on the bodies $\gamma_\ell$ by numerically computing the shear stress as described in Section~\ref{sec:shearStress}.

\subsection{Pressure}
\label{sec:pressure}
To compute the pressure at $\xx \in \Omega$, we follow the same procedure used for the shear stress and vorticity.  We first compute the pressure of the double-layer potential for $\xx \in \Omega$ and add the appropriate jump if $\xx \in \Gamma$. The pressure of the double-layer potential~\eqref{eqn:velocityDLP} at $\xx \in \Omega$ is
\begin{align}
  p^{\DD}(\xx) = -\frac{1}{\pi}\int_{\bd\Omega} \frac{1}{\rho^2}
    \left(I - 2 \frac{\rr \otimes \rr}{\rho^2}\right) 
    \nn \cdot \eeta(\yy) \,ds_\yy \, .
    \label{eqn:pressureDLP}
\end{align}
To compute~\eqref{eqn:pressureDLP} numerically, we apply the trapezoid rule with an upsampling scheme described in Section~\ref{sec:NSI}.  The pressure due to the Stokeslets is
\begin{align}
  p^S(\xx) = \sum_{\ell=1}^{M}\frac{\rr \cdot
  \llambda_\ell}{2\pi\rho^2}\, ,
\end{align}
and the pressure due to the rotlets is zero.  Therefore, the pressure is
\begin{align}
  p(\xx) = p^{\DD}(\xx) + p^{S}(\xx)\, , \quad \xx \in \Omega \, .
\end{align}

For target points on the boundary of the geometry, $\xx \in \bd \Omega$, there is an additional term resulting from a jump in the pressure double-layer potential~\cite{poz1992}
\begin{align}
  p(\xx)= \pderiv{\eeta}{\ss} \cdot \ss + p^{\DD}(\xx) + 
              p^{S}(\xx)\, , \quad \xx \in \bd\Omega \, .
  \label{eqn:pressureJump}
\end{align}
Similar to the shear stress, the singularity of the kernel~\eqref{eqn:pressureDLP} is too strong when $\xx \in \bd\Omega$ for the trapezoid rule to be directly applied.  Therefore, we again subtract a constant density function and write the pressure as
\begin{align}
  p^\DD(\xx) = -\frac{1}{\pi}\int_{\bd\Omega} \frac{1}{\rho^2}
    \left(I - 2 \frac{\rr \otimes \rr}{\rho^2}\right) 
    \nn \cdot (\eeta(\yy) - \eeta(\xx)) \,ds_\yy\, , 
    \quad \xx \in \bd\Omega\, .
  \label{eqn:pressureOddEven}
\end{align}
Then, identical to how the shear stress is computed in~\eqref{eqn:stressOddEven}, odd-even integration is used to evaluate~\eqref{eqn:pressureOddEven}, and the jump condition in~\eqref{eqn:pressureJump} is computed with Fourier differentiation.  

\subsection{Drag}
\label{sec:drag}
Once the pressure and shear stress are computed, the drag on a single body with boundary $\gamma_\ell$ is
\begin{align}
\FFD = \int_{\gamma_\ell} p \nn + \tau \ss \,ds\, .
\label{eqn:drag}
\end{align}
The reader may notice an unusual sign convention for pressure in this formula, which is simply due to the normal direction pointing {\em into} the bodies.

\subsection{Near-singular integration}
\label{sec:NSI}
The layer potentials for the velocity~\eqref{eqn:velocityDLP}, vorticity~\eqref{eqn:vorticityDLP}, and pressure~\eqref{eqn:pressureDLP} need to be evaluated at points $\xx \in \Omega$. We write the general layer potential as
\begin{align}
  \phi(\xx) = \int_{\bd\Omega} K(\xx,\yy) \eeta(\yy) \, ds_\yy\,, \quad \xx \in \Omega\, .
  \label{eqn:genericLP}
\end{align}
Since the integrand is periodic and smooth, the trapezoid rule guarantees spectral accuracy.  However, when the target point $\xx$ is close to $\bd\Omega$, the integrand is nearly-singular, and the accuracy of the trapezoid rule deteriorates unless a large number of discretization points are used~\cite{bar2014}.  We are only interested in computing the velocity, vorticity, and pressure to visualize the flow, so it is not necessary that we resolve the integrands for arbitrarily close target points.

The layer potential~\eqref{eqn:genericLP} is approximated with the trapezoid rule 
\begin{align}
  \phi(\xx) = \sum_{j=1}^N q_j K(\xx,\yy_j)\, ,
  \label{eqn:layerPotential}
\end{align}
where $q$ contains the quadrature weights, Jacobian, and density function, and $N$ depends on the distance between $\xx$ and $\bd\Omega$.  Based on numerical experiments, the error of the $N$-point trapezoid rule is acceptable for points $\xx$ that are more than five arclength spacings from $\bd\Omega$. For points that are closer than five arclength spacings to $\bd\Omega$, we use an upsampled trapezoid rule.  The required upsampling of the geometry and the density function is done with a Fourier interpolant.  To minimize the number of discretization points that are used to approximate the layer potential~\eqref{eqn:layerPotential}, we use the same upsampling rate for all target points in dyadic intervals (Table~\ref{tbl:upsampling}).  For points that are closer than $5/16$ of an arclength from $\bd\Omega$, a value of 0 is simply assigned to $\phi$.
\newcolumntype{V}{>{\centering\arraybackslash}m{4em}}
\newcolumntype{N}{@{}m{0pt}@{}}
\begin{table}[htpb]
\begin{center}
\caption{\label{tbl:upsampling}The upsampling rate of our near-singular integration scheme.  $d(\xx,\bd\Omega)$ is the distance between $\xx \in \Omega$ and $\bd\Omega$, and $h$ is an arclength spacing.}
\vspace{0.3 pc}
\label{table:arangle}
\begin{tabular}{|V|V|V|V|V|V|N}
  \hline
  $d(\xx,\bd\Omega)$ &
  $\left[5h,\infty\right]$ &
  $\left[\frac{5}{2}h,5h\right]$ &
  $\left[\frac{5}{4}h,\frac{5}{2}h\right]$ & 
  $\left[\frac{5}{8}h,\frac{5}{4}h\right]$ &
  $\left[\frac{5}{16}h,\frac{5}{8}h\right]$ & \\ [2ex] 
  \hline
  Upsample Rate & 1 & 2 & 4 & 8 & 16 & \\
  \hline
\end{tabular}
\end{center}
\end{table}

If eroding bodies are also sedimenting, bodies will come too close to one another for upsampling to be sufficient.  Therefore, future work will apply one of several near-singular integration schemes available in the literature~\cite{qua-bir2014a, klo-bar-gre-one2013, bar-wu-vee2015, hel-oja2008a}.

\section{Results: Single-body erosion}
\label{s:SingleResults}
With the description of the methods complete, we now present numerical results on bodies eroding in Stokes flow. We first discuss the case of a single body, as certain features can be predicted analytically and thus used for validation. In Section~\ref{s:MultiResults} we will discuss the erosion of multiple bodies and highlight important differences that arise compared to the single-body case.

\subsection{A single body eroding in Stokes flow}

\begin{figure}
\begin{center}
\includegraphics[width = 0.80 \textwidth]{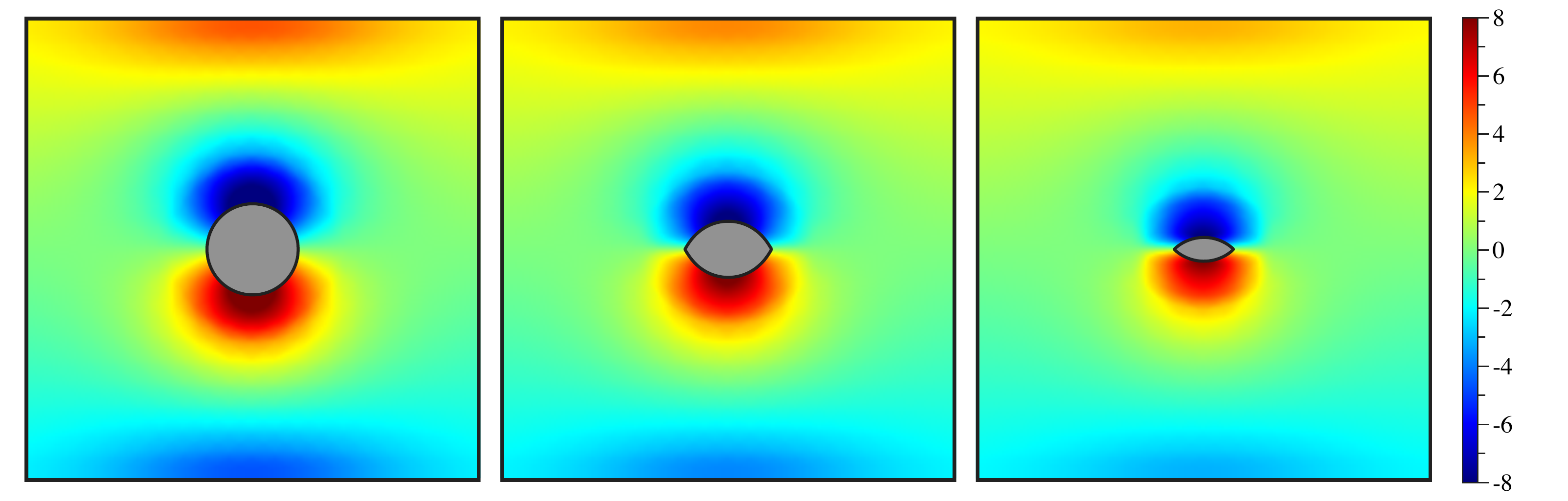}
\caption{A single body eroding in Stokes flow. As the body erodes, it develops a corner at its front and rear stagnation points, while remaining for-aft symmetric due to time-reversibility. Color indicates the vorticity of the surrounding flow, which, when evaluated on the boundary, gives the local shear stress. In this simulation, we use $N_\iin = 1024$ discretization points on the body, $N_\out = 1024$ points on the outer boundary, a time-step of $\Dt = 10^{-6}$, and smoothing parameters of $\eps = 10/1024$ and $\sigma = 10/1024$.}
\label{01bodseq} 
\end{center}
\end{figure}

Figure~\ref{01bodseq} shows the erosion of a single body, numerically simulated with $N_\iin = 1024$ discretization points and shown at three evenly spaced times. Here, color represents the vorticity of the surrounding flow. When evaluated on the boundary, vorticity corresponds to the shear stress and thus the local erosion rate. As seen in the left-most snapshot, the body is initially circular and has the strongest vorticity at its top and bottom. The flow is for-aft symmetric as a result of the time-reversibility of the Stokes equations, with vorticity vanishing at the body's front and rear stagnation points. At the next instance, the body has lost material due to erosion, but remains for-aft symmetric due to the flow symmetry. The greatest material loss has occurred at the body's top and bottom, where shear is strongest. Interestingly, the body has developed well-defined corners at its front and rear stagnation points. This geometric feature is a consequence of the erosion law, which depends on the {\em absolute} shear stress: wherever the stress vanishes, $\atau$ exhibits a corner, which manifests as a corner in the geometry. As shown in the third snapshot, this corner persists as the bodies becomes even more slender and continues to shrink.  

\begin{figure}
\begin{center}
\includegraphics[width = 0.8 \textwidth]{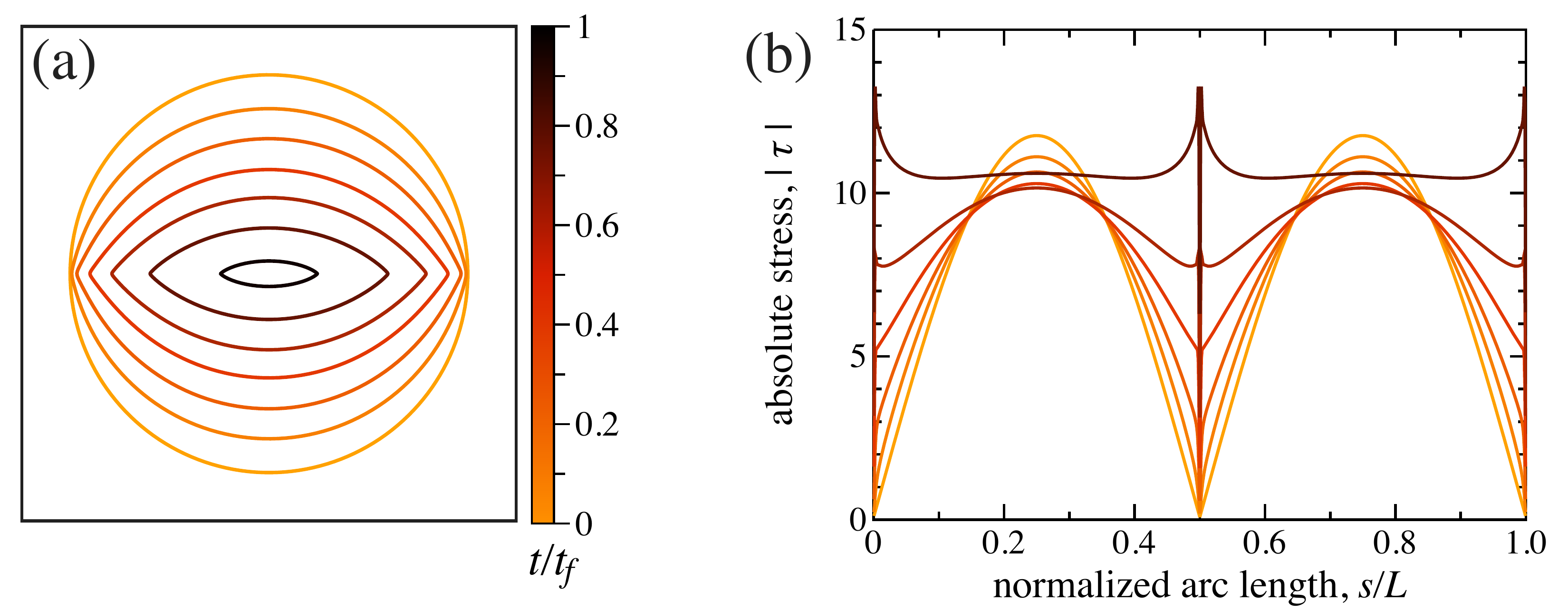}
\caption{
Shape and stress evolution: (a) Interfaces of the eroding body from Fig.~\ref{01bodseq} shown at evenly-spaced time intervals. Color indicates normalized time, $t/t_f$, where $t_f = 1.79 \times 10^{-2}$ is the time at which the body completely vanishes. (b) The distribution of absolute shear stress, $\atau$, along each of these interfaces. Normalized arclength, $s/L$, runs counterclockwise from the rear stagnation point. The stress initially varies over the body's surface, but becomes more uniform over time.}
\label{shrink_intface}
\end{center}
\end{figure}

To better understand the shape evolution, Fig.~\ref{shrink_intface}(a) shows several superimposed interfaces at evenly spaced times. The interfaces are color-coded by normalized time $t/t_f$, where $t_f = 1.79 \times 10^{-2}$ is the time at which the body vanishes in the simulation. As seen in Fig.~\ref{shrink_intface}(a), the body converges to a fairly slender, for-aft symmetric shape that has corners at its front and back---somewhat reminiscent of an American football.  Fig.~\ref{shrink_intface}(b) shows how the absolute shear stress, $\atau$, is distributed along each of these interfaces. The normalized arclength $s/L$ runs counterclockwise around the body, with $s/L = $ 0 and 0.5 corresponding to the rear and front stagnation points respectively. At early times (yellow curves), the stress varies significantly, with more stress concentrated at the top and bottom of the body. However, as erosion reshapes the body, the stress becomes more uniformly distributed over the entire surface (red curves). The absolute stress $\atau$ shows slight peaks near the stagnation points as the body nears the final stages of erosion (dark red). Closer inspection reveals that the stress actually vanishes {\em at} the stagnation points, but in a region that becomes increasingly narrow with time.

\subsection{Scaling laws for area and drag}
\label{sec:scaling}

Figure~\ref{shrink_intface}(a) shows that the spacing between successive interfaces becomes larger with time, suggesting that the overall shear stress grows as the body vanishes. This observation is consistent with the well-known 3D Stokes drag law, in which stress, $\tau \sim \mu \umax / a$, scales inversely proportional to body-scale $a$ (e.g.~for a sphere, $a$ is the radius). This law must be modified in two dimensions, though. The modification is made non-trivial by subtleties of the 2D Stokes paradox, but can be correctly deduced by noting that, as the body vanishes, its drag must be consistent with slender-body theory~\cite{batchelor1970slender, MooreJFM2012}.  For a slender body in flow perpendicular to its long-axis $\ell$, the drag is proportional to $\mu \umax \ell/ \log(\ell/a)$, where $a \ll \ell $ is a length scale of the cross-section. This law implies that the stress-per-unit-span on a 2D cross-section is approximately
\begin{align}
  \label{eqn:stresslaw}
  \tau \sim \frac{\mu \umax}{a \log(\ell/a)}\, , \qquad \text{for } a \ll \ell \, .
\end{align}
In interpreting our 2D simulations, we will treat the span, $\ell$, as fixed. Recall that $L$ denotes the total perimeter of the 2D cross-section. Since $L$ and $a$ each represent a characteristic scale of the cross-section, they can often be interchanged in scaling analysis.  

\begin{figure}
\begin{center}
\includegraphics[width = 0.9 \textwidth]{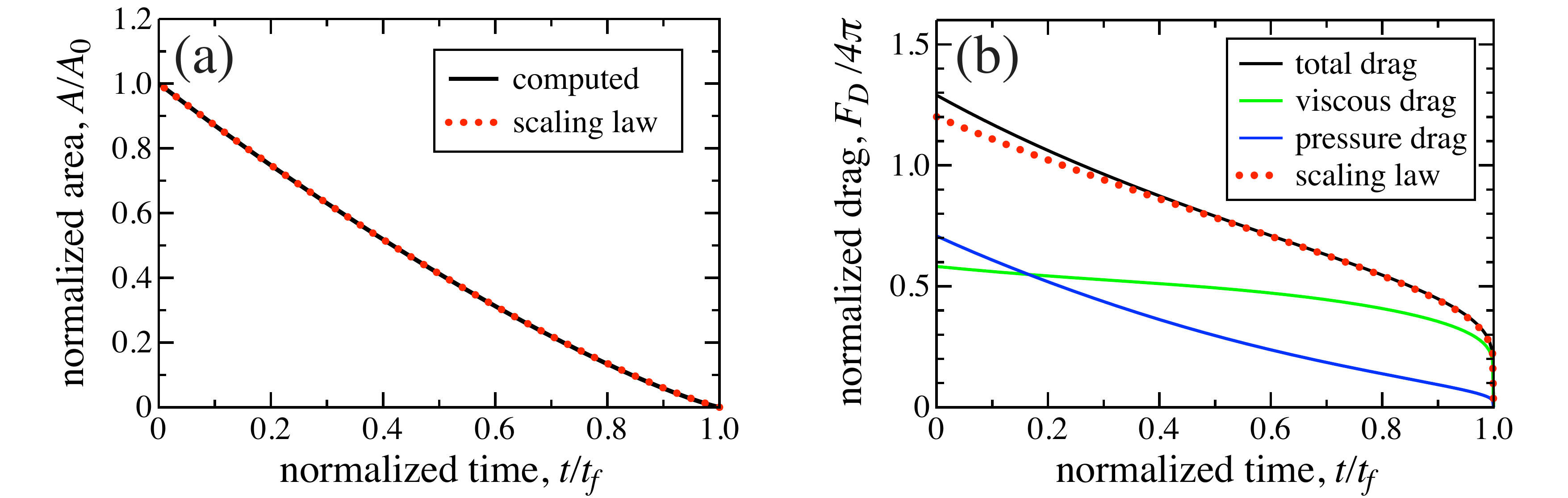}
\caption{Area and drag: (a) The area of the eroding body as measured in the simulations (black) shown against scaling law~\eqref{area_predict} (dotted red). The prediction and measurements agree closely. (b) The total horizontal drag exerted by body (black), along with the viscous (green) and pressure (blue) contributions. Slender-body theory yields the drag estimate~\eqref{dragscaling} (dotted red), which agrees well with the simulation as $t \to t_f$.}
\label{area_drag}
\end{center}
\end{figure}

The body's rate of area reduction is given by integrating the absolute shear stress around the boundary
\begin{align}
\dot{A} = - \int_{\gamma} \atau \, ds \,  = - \mean{\atau} L \, .
\end{align}
Thus, the rate of area reduction is given by the product of the mean shear stress $\mean{\atau}$ and the perimeter $L$, which is an {\em exact} relation. Inserting scaling law~\eqref{eqn:stresslaw} gives the estimate
\begin{align}
\dot{A} \sim \frac{ - \mu \umax}{\log(\ell/a)} \qquad \text{as } a \to 0 \, . 
\end{align}
Noting that $A \sim a^2$ allows one to solve the ODE, giving
\begin{align}
\label{area_predict}
\frac{A}{A_0} \left( 1 - c_1 \log{\frac{A}{A_0}} \right) \sim \frac{t_f - t}{t_f} \quad \text{as } t \to t_f \, .
\end{align}
Here, $t_f$ is the time of vanishing, $A_0$ is the body's initial area, and $c_1$ is a free constant. The above is an {\em implicit} relationship for how area scales down in time, and we can use this prediction to test our numerical simulation. In Fig.~\ref{area_drag}(a) we compare the area measured in our numerical simulation (black curve) against the above scaling law (red dotted curve) with $c_1 = 0.24$ a fit parameter. The two agree remarkably well throughout the entire simulation, thus co-validating the simulation and asymptotic results.

As the body shrinks in accordance with~\eqref{area_predict}, the drag that it exerts decreases at a related rate, offering a second check on our numerics. In Fig.~\ref{area_drag}(b), we show the horizontal drag, $F_D$, on the simulated body as it changes in time (black curve). The drag is calculated using~\eqref{eqn:drag} and, in the context of slender-body theory, this should be interpreted as the drag-per-unit-span. As seen in the figure, $F_D$ decreases gradually at first, then abruptly drops to zero in the final moments of the body's existence. Though it may seem peculiar, the sudden drop is consistent with the 2D Stokes law in which drag vanishes {\em logarithmically} with body size. Additionally,~\eqref{eqn:drag} permits us to calculate the separate contributions of pressure and viscous drag, which we show by the blue and green curves in Fig.~\ref{area_drag}(b) respectively. Initially, pressure and viscous drag are comparable, with the circular geometry having a slightly higher level of pressure drag. As the body erodes and assumes a more slender form, though, the pressure drag decreases more rapidly, leaving viscous drag to dominate. This behavior is associated with convergence to a drag-minimizing profile, which is discussed in Section~\ref{LimitingShape}.

As a check on these results, we can use slender-body theory to estimate the drag-per-unit-span
\begin{align}
\label{dragscaling}
\frac{F_D}{4 \pi} \sim \frac{ \mu U}{\log(\ell/a)}\, ,	\qquad \text{for } a \ll \ell \, .
\end{align}
For a {\em circular} cylinder, the normalization factor $4 \pi$ produces an {\em exact} relationship (within the context of slender-body theory), where $a$ is taken to be the radius. For other shapes, the formula will be approximate and we define $a = \sqrt{A/\pi}$ as the characteristic length scale since it generalizes the case of a circle. In Fig.~\ref{area_drag}(b) we show this prediction (red dotted curve), where $\ell = 2.3 a_0$ is a fit parameter and $a_0 = 0.2$ is the radius of the initially circular body. The scaling law and measured drag agree closely as $t \to t_f$, confirming that the simulation accurately captures the slender limit as the body disappears.

\subsection{Local analysis for the limiting shape}
\label{LimitingShape}
 
\begin{figure}
\begin{center}
\includegraphics[width = 0.4 \textwidth]{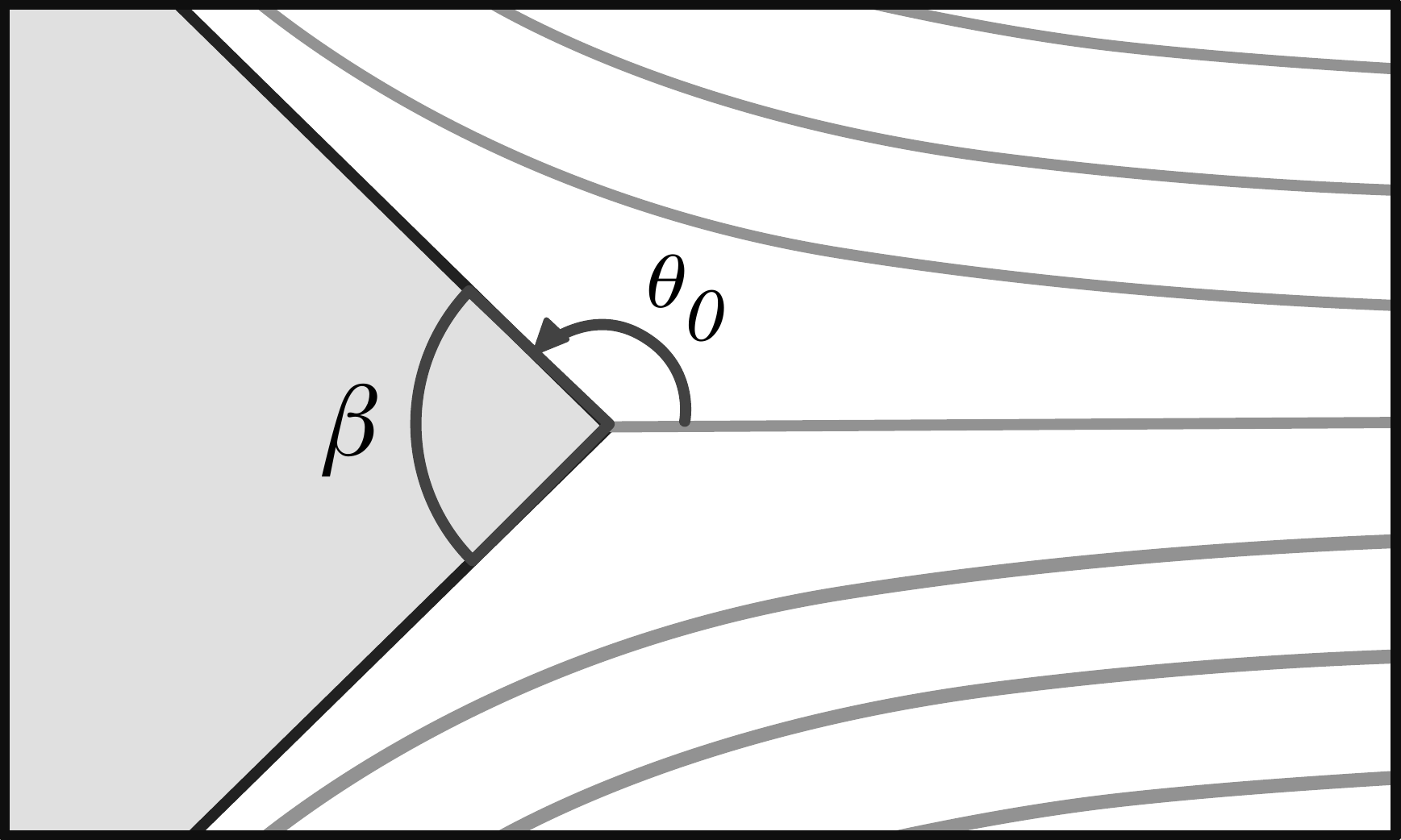}
\caption{Schematic of 2D Stokes flow near a corner: The flow around an infinite wedge of opening angle $\oangle$ is described by a stream-function in polar coordinates $\psi(r,\theta)$, for $\theta \in [-\theta_0, \theta_0]$.}
\label{corner}
\end{center}
\end{figure}
 
With the vanishing rate understood, we now seek to quantify the limiting shape observed in Fig.~\ref{shrink_intface}(a). In particular, consider the corner formed at the front and rear of the body. Locally, we can neglect body curvature and consider the problem of Stokes flow around an infinite wedge; see Fig.~\ref{corner} for a schematic. We seek to predict the wedge's opening angle, $\oangle$, by appealing to the principle that, in long time, erosion uniformly distributes the shear stress along the surface~\cite{moo-ris-chi-zha-she2013, moore2017riemann}. For a wedge of arbitrary angle, the surrounding flow can be determined exactly in polar coordinates $(r, \theta)$ where $\theta \in [-\thb,\thb]$ and $\thb = \pi - \oangle/2$ (in this section alone, $\theta$ denotes the polar coordinate and not the local tangent angle). 

The corresponding stream function is biharmonic, $\Delta^2 \psi = 0$, and the ansatz $\psi = r^{\lambda}f(\theta)$ gives the 4th-order ODE~\cite{poz1997}
\begin{align}
\label{fODE}
f'''' + 2(\lambda^2 - 2 \lambda + 2)f'' + \lambda^2(\lambda-2)^2 f = 0 \, .
\end{align}
This ODE is subject to no-slip boundary conditions along the wedge $f(\thb) = f'(\thb) = 0$ and odd symmetry about the horizontal axis, $f(0) = f''(0) = 0$. For arbitrary opening angle, ODE~\eqref{fODE} and the boundary conditions comprise an eigenvalue problem for $\lambda$, whose solution furnishes the stream function and the corresponding shear stress on the wedge surface $\tau = \mu r^{\lambda-2} f''(\thb)$. We, however, seek to determine the particular opening angle that produces {\em uniform} shear stress and therefore set $\lambda = 2$.  The odd-symmetric solution to~\eqref{fODE} corresponding to $\lambda=2$ is
\begin{align}
  f(\theta) = B \sin (2 \theta) + C \theta \, .
\end{align}
Imposing no slip along the wedge surface gives the condition
\begin{align}
  \sin(2 \thb) - 2 \thb \cos(2 \thb) = 0 \, .
\end{align}
This nonlinear equation has solution $\thb \approx 129^\circ$, which implies that the uniform-stress wedge has an opening angle of $\oangle \approx 102^\circ$.

\subsection{Fixed-area simulations}

\begin{figure}
\begin{center}
\includegraphics[width = 0.85 \textwidth]{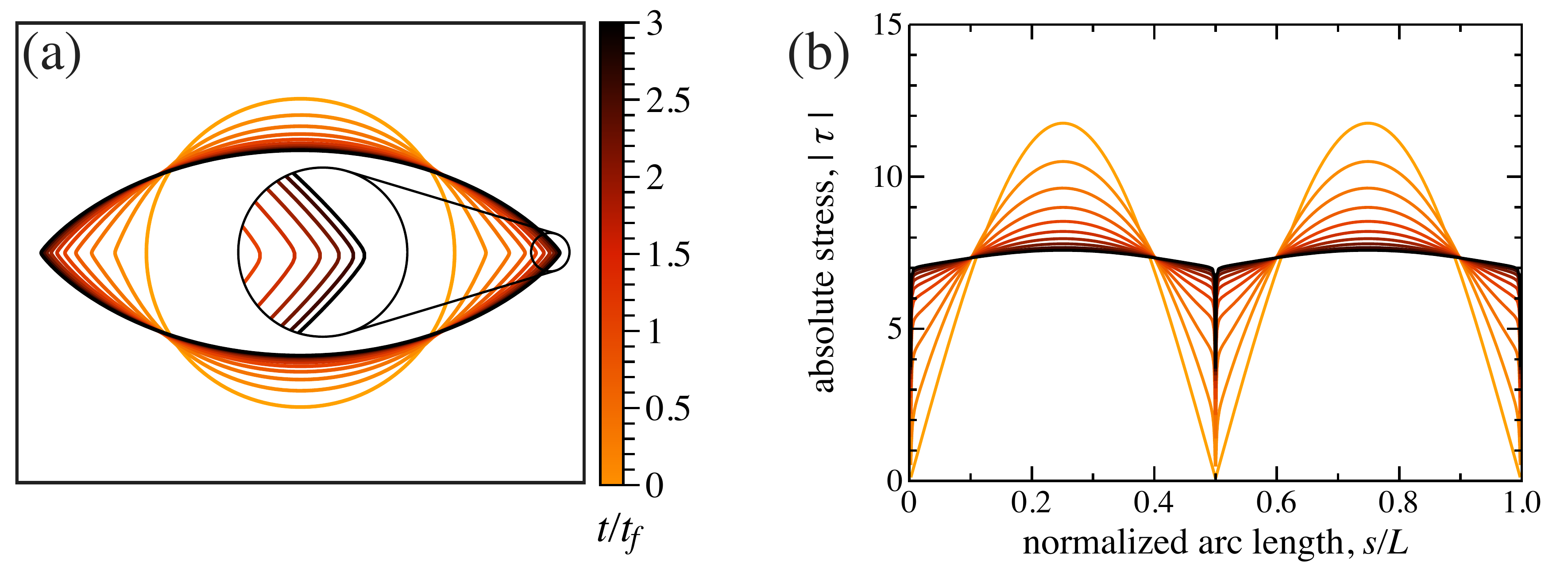}
\caption{Fixed-area simulations: (a) Interfaces of a body eroding under shear stress but with the area artificially kept constant. Inset: zoom near the rear stagnation point, showing that the corner is actually smooth on a fine scale. (b) The shear-stress distribution on each of these interfaces. As the body erodes and converges to a terminal form, the stress becomes more uniformly distributed over the entire surface. Here time is normalize by the same value $t_f = 1.79 \times 10^{-2}$ from Fig.~\ref{shrink_intface}.}
\label{fixed_intface}
\end{center}
\end{figure}

We now aim to measure the opening angles formed in the numerical simulations to compare against the above prediction. However, Section~\ref{sec:scaling} shows that the magnitude of the shear stress diverges as the body vanishes, which causes our simulations to lose accuracy in the final, critical moments when the body is expected to approach its limiting shape. To circumvent this problem, we remove the effect of shrinking by forcing $\Vn$ to have zero mean,
\begin{align}
\Vn \leftarrow \Vn - \mean{\Vn} \, .
\end{align}
This modification allows erosion to sculpt the shape as it would normally, but with the area artificially held fixed.

We show in Fig.~\ref{fixed_intface}(a) several interfaces of this modified erosion process, with the same spatial resolution, $N_\iin = 1024$, as before. Notice the body obtains a similar profile as observed previously, but now with fixed area. Importantly, more numerical accuracy is available to quantify the terminal state. The inset in Fig.~\ref{fixed_intface}(a) zooms on the rear corner and demonstrates that, due to our numerical regularization, the body is actually smooth on a fine scale. We show in Fig.~\ref{fixed_intface}(b) the stress distributions on the same interfaces. As observed previously, the stress becomes more evenly distributed as erosion reshapes the body. The slight rise in shear near the front/rear corners that was observed in Fig.~\ref{shrink_intface}(b) is no longer present, suggesting that this effect was due to loss of accuracy. Rather, the late-time stress distributions in Fig.~\ref{fixed_intface}(b) decrease monotonically to zero near the stagnation points. In this figure, time is normalized by the vanishing time $t_f = 1.79 \times 10^{-2}$ from the shrinking simulations (i.e.~Fig.~\ref{shrink_intface}). Curiously, about 3 multiples of $t_f$ are required to reach a visible steady state, which contrasts with previous studies on high-Reynolds-number erosion in which a steady state usually emerged before the body disappeared~\cite{moo-ris-chi-zha-she2013}.  
\begin{figure}
\begin{center}
\includegraphics[width = 0.85 \textwidth]{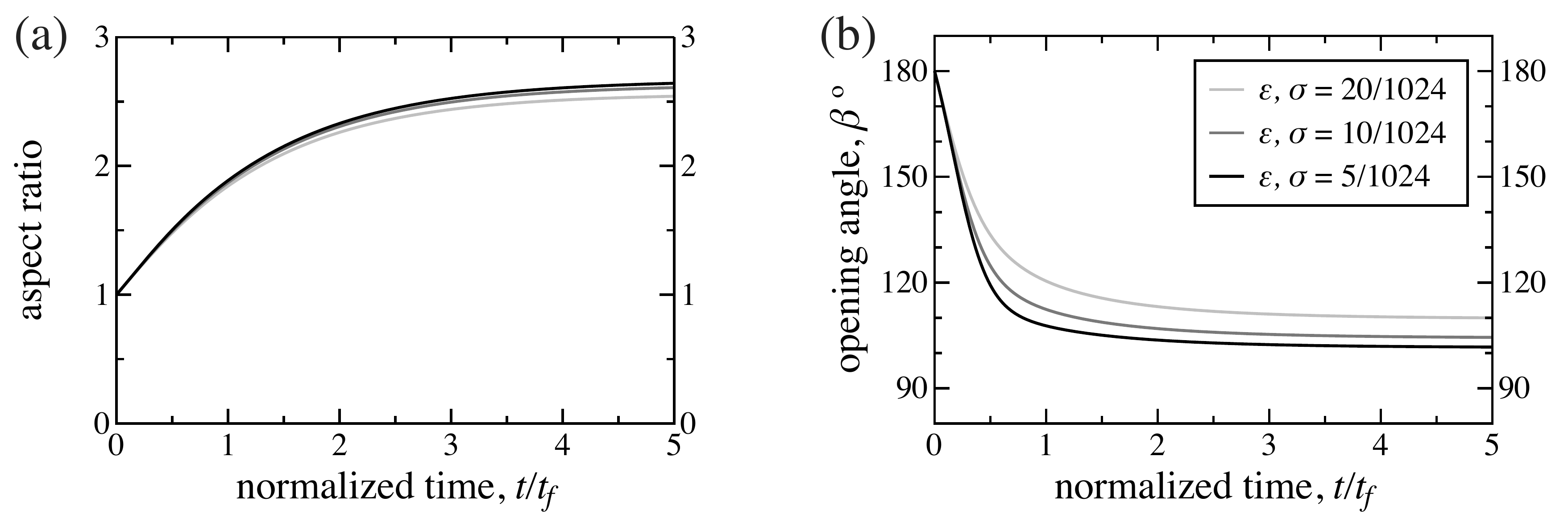}
\caption{Shape characteristics of an eroding body: (a) The aspect ratio (width to height) versus time. (b) The measured opening angle versus time. For both measurements, we show three cases of the smoothing parameters: $\eps = \sigma = 20/1024, 10/1024,$ and 5/1024. In the finest simulation ($5/1024$) the aspect ratio converges to 2.7 and the opening angle converges to the predicted value of $102^{\circ}$ in long time.}
\label{fig:arangle}
\end{center}
\end{figure}

These fixed-area simulations enable us to precisely characterize the changing geometry as it approaches the terminal state. In particular, we measure the aspect ratio of the body (the width to height) and the opening angle of the front/rear corners. In these tests, we keep the spatial resolution fixed, $N_\iin = 1024$, and show three different cases of the smoothing parameters, $\eps$ and $\sigma$, to assess how the numerical regularization influences shape. In the three cases, we keep $\eps = \sigma$ always and set them to the values 20/1024, 10/1024, and 5/1024. Loosely, these values correspond to smoothing over the nearest 20, 10, and 5 grid points, respectively. Figure~\ref{fig:arangle}(a) shows the aspect ratio as it varies in time. The initially circular geometry has an aspect ratio of unity, and, as erosion carves a more slender form, the aspect ratio increases monotonically. The three simulations are all very close to one another, demonstrating convergence as $\eps, \sigma \to 0$. In the simulation with the least smoothing ($\eps = \sigma = 5/1024$), the aspect ratio tends to a value of 2.7 in long time.

Next, we aim to measure the opening angle of the front/rear corner. Since these simulated bodies are actually smooth on a fine scale, we must have a procedure to estimate the apparent angle that results from the high (but finite) curvature at the stagnation points. To this end, we fit the tangent angle $\theta(\alpha)$ with a 7th-degree polynomial over a region excluding the stagnation points ($\alpha = 0$ and 0.5). We then extrapolate the polynomial to the stagnation points, which gives an estimate of the angle there. Although this procedure is subject to some uncertainty (to be examined below), we have found it to give robust angle measurements.  Figure~\ref{fig:arangle}(b) shows the opening angle, $\beta$, as it varies in time for the three cases of smoothing parameters.  For the initially circular geometry, our procedure yields the correct opening angle of $180^{\circ}$, i.e.~a flat front face.  As time proceeds, erosion carves a sharper corner and $\beta$ decreases monotonically. Here, the three cases of smoothing parameters are more easily distinguished, implying that the opening angle is more sensitive to the level of numerical smoothing than the aspect ratio is. The curves nevertheless show convergence as $\eps, \sigma \to 0$. In the finest simulation ($\eps = \sigma = 5/1024$), $\beta$ converges to the value $102^{\circ}$, which agrees precisely with the analytical prediction from Section~\ref{LimitingShape}. 

To more thoroughly examine the influence of the smoothing parameters, we report in Table~\ref{table:arangle} the long-time values of the aspect ratio and opening angle for each of the three cases. The table shows that both the aspect ratio and opening angle converge at a rate that is nearly linear as $\eps, \sigma \to 0$. In addition, we show an estimate of the uncertainty involved in our angle measurement, obtained by varying the degree of the fitting polynomial and the size of the excluded region around the stagnation points. The estimates show that as $\eps, \sigma \to 0$ the angle-measurements become more precise and converge to the predicted value of $102^{\circ}$.  

\begin{table}
\begin{center}
\caption{The long-time aspect ratio and opening-angle measurements for three cases of numerical smoothing. In these runs, we fix the resolution $N_\iin = 1024$, and vary the two smoothing parameters, $\eps, \sigma$. As $\eps, \sigma \to 0$, the aspect ratio tends to roughly 2.7 and the opening angle converges to the predicted value of $102^{\circ}$. Furthermore, the uncertainty in the angle measurement becomes smaller.
} 
\vspace{0.3 pc}
\begin{tabular}{c c c}
\hline
\hspace{0.5pc} smoothing parameters: $\eps$ and $\sigma$
\hspace{0.5pc} & final aspect ratio 
\hspace{0.5pc} & final opening angle \\
\hline
20/1024		& 2.55	& $110^\circ \pm 6^\circ$	\\
10/1024		& 2.62	& $104^\circ \pm 4^\circ$	\\
5/1024		& 2.65	& $102^\circ \pm 2^\circ$	\\
\hline
\end{tabular}
\end{center}
\end{table}

In summary, we have found that a single-body eroding in 2D Stokes flow converges to a terminal form characterized by nearly uniform shear stress along its surface. The final shape is fairly slender, with an aspect ratio of 2.7, and distinguished by  corners of opening angle $102^{\circ}$ at the front and rear. This shape is closely related to the drag-minimizing profile found by Pironneau (1973)---the geometry that minimizes drag in 3D Stokes flow subject to the constraint of fixed volume~\cite{pir1973}. Using variational principles, Pironneau showed that the drag-minimizing form has constant shear (or vorticity) along its surface. The resulting {\em axisymmetric} shape has an aspect ratio of 2.1 and front/rear opening-angles of $120^{\circ}$~\cite{pir1973, mit-spa2016}.  Interestingly, 3D simulations of viscous erosion showed close resemblance between the shape formed by a single eroding body and this Pironneau form~\cite{mit-spa2016}. It was found that the eroding body does not converge {\em precisely} to the Pironneau shape because the long-time surface-stress distribution created by erosion is only {\em approximately} uniform. Rather than strict uniformity, variations in body curvature must be offset by slight variations in $\Vn$ for self-similar reduction~\cite{pir1973, mit-spa2016}. Pironneau's argument that the drag-minimizing form has constant surface stress applies to 2D as well. Therefore, the shape observed in our simulation is approximately the form that minimizes drag subject to fixed area, i.e.~the 2D counterpart to Pironneau's shape. How close this shape is to the true drag-minimizer is related to how close the stress distributions are to uniform, as can be seen in Fig.~\ref{fixed_intface}(b).  

\subsection{Convergence test}
As a final validation, we now verify that our method converges with the expected second-order accuracy in time. For this, we run the same single-body simulation shown in Figs.~\ref{01bodseq} and~\ref{shrink_intface} (except with a coarser $N_\iin=256$) to a stopping time of $t_s = 10^{-2}$. This stopping time is 56\% of the vanishing time, and was selected so that the body will have eroded substantially but is not small enough to compromise accuracy (due to the stress diverging).  We assess error by measuring the $L^2$-difference in the body shape obtained from simulations of successive resolution, with $\Delta t$ halved at each stage. In Table~\ref{convtab}, we report these errors along with the order of accuracy, obtained by Richardson-extrapolation-type calculations. The table confirms that our method converges with expected second-order accuracy. The last column of the table shows the CPU time taken for the simulation to run on a single-processor machine.

\begin{table}
\begin{center}
\caption{Convergence test: We run single-body erosion simulations to a stopping time of $t_s = 10^{-2}$ (56\% of the vanishing time) and measure the $L^2$-difference in shapes obtained with successive temporal resolutions. The tests confirm the method to be second-order accurate in time. The last column shows the measured CPU time which increases linearly as the time-step is refined.
}
\vspace{0.3 pc}
\label{convtab}
\begin{tabular}{l l l l}
\hline
\hspace{0.0pc} $\Delta t/t_s$
\hspace{0.5pc} & error 
\hspace{0.5pc} & order
\hspace{0.5pc} & CPU time \\
\hline
1/50		& 1.68E-5		& --		& 53 secs     	\\
1/100	& 4.84E-6		& 1.79	& 1.8 mins   	\\
1/200	& 1.30E-6		& 1.90	& 3.5 mins  	\\
1/400	& 3.36E-7		& 1.95	& 6.9 mins  	\\
1/800	& 8.54E-8		& 1.98	& 14 mins   	\\
1/1600	& 2.15E-8		& 1.99	& 28 mins  	\\
1/3200	& 5.40E-9		& 1.99	& 56 mins    	\\
1/6400	& --			& --		& 111 mins	\\
\hline
\end{tabular}
\end{center}
\end{table}

\section{Results: Multibody erosion}
\label{s:MultiResults}

\subsection{Fixed pressure-drop condition}

With the numerical method validated, we are now ready to simulate the erosion of multiple bodies. We will first make one conceptual change to the model. For a realistic porous medium, the condition of a fixed incoming flow rate, $\umax$, is somewhat artificial. In groundwater flow, for example, it is more common to encounter an upstream (recharge) region of given hydraulic head that connects to a downstream (discharge) region of another head. Motivated by this example, we will impose the condition that the pressure drop across the computational domain remains fixed in time. Consequently, we would expect the flow rate, $\umax$, to increase over time as the bodies erode inside the cell and thus provide less resistance.

We will impose the fixed-pressure-drop condition by dynamically controlling the incoming flow rate, $\umax$. We first compute the Stokes flow with $\umax = 1$ prescribed and measure the corresponding pressure at an upstream ($x = -2$) and downstream ($x=2$) location. We next calculate the vertically-averaged upstream and downstream pressures, defined as
\begin{align}
p_{-} = \frac{1}{2} \int_{-1}^{1} p(-2,y) \, dy \, , \qquad
p_{+} = \frac{1}{2} \int_{-1}^{1} p(+2,y) \, dy \, .
\end{align}
We then rescale the computed density function, Stokeslets, and rotlets  by the pressure difference
\begin{align}
\eeta \leftarrow \frac{\eeta}{(p_{-} - p_{+})}\, , \quad 
\llambda \leftarrow \frac{\llambda}{(p_{-} - p_{+})}\, , \quad 
\xi \leftarrow \frac{\xi}{(p_{-} - p_{+})}\, .
\end{align}
The velocity and pressure fields depend linearly on $\eeta$, $\llambda$, and $\xi$, so this rescaling maintains a {\em fixed} pressure difference across the subdomain $-2 \le x \le 2$. Since the Stokes equations are linear, the rescaled flow remains a solution to the PDE system (valid to the same order of accuracy). Furthermore, since the Stokes equations are time-independent, rescaling $U$ by a different value at each time step does not introduce any additional error.

\subsection{Validation}

To validate the multibody simulations with this new fixed-pressure-drop condition, we report a convergence test in Table~\ref{convtab2}. In this test, we simulate the erosion of 20 bodies to a stopping time of $t_s = 4 \times 10^{-3}$. This stopping time is 24\% of the time required for {\em all} bodies to vanish, and, in fact, a few of the bodies will have already vanished by $t_s$. To assess error, we measure the aggregate $L^2$-difference in the shape of {\em all} of the bodies, comparing simulations of successive temporal resolution. Note that this metric is sensitive to numerical errors in any of the simulated bodies. The table confirms that our method retains second-order accuracy in time, thus validating our rescaling procedure and confirming that multibody interactions are handled correctly by the Stokes solver. The last column of the table shows the computational time of each simulation, demonstrating a linear increase with refinement of the time step. Naturally, these 20-body simulations require significantly more CPU time than do the single-body simulations.

\begin{table}
\begin{center}
\caption{Convergence test for multiple-body erosion with a fixed-pressure-drop condition. In this test, we use 20 bodies with $N_{\iin} = 128$ grid points on each. We run the simulations to a stopping time of $t_s = 4 \times 10^{-3}$ (24\% of the time it takes for {\em all} bodies to vanish). We measure the {\em combined} $L^2$-difference in the shapes of {\em all} bodies. The table confirms that the method maintains second-order accuracy in time, even with multiple bodies and a time-varying inflow rate $\umax$.
}
\vspace{0.3 pc}
\label{convtab2}
\begin{tabular}{l l l l}
\hline
\hspace{0.0pc} $\Delta t/t_s$
\hspace{0.5pc} & error 
\hspace{0.5pc} & order
\hspace{0.5pc} & CPU time \\
\hline
1/10     	& 2.67E-5  	& --        	& 14 mins  	\\
1/20     	& 7.23E-6  	& 1.88 	& 27 mins  	\\
1/40     	& 1.95E-6  	& 1.89 	& 52 mins  	\\
1/80     	& 4.90E-7  	& 1.99 	& 1.7 hours	\\
1/160     	& 1.25E-7  	& 1.97  	& 3.4 hours	\\
\hline
\end{tabular}
\end{center}
\end{table}

\begin{table}
\begin{center}
\caption{The number of GMRES iterations required to reach a fixed tolerance of $10^{-8}$ in our simulations. The left side shows a case of fixed geometry, 10 circular bodies, with increasing spatial resolution, $N_\iin$, demonstrating that the iteration count is independent of resolution. In the right side, we fix the resolution and sequentially increase the number of bodies (all circular) from 10 to 50. The iteration count increases roughly linearly with the number of bodies.  }
\vspace{0.3 pc}
\label{itertab}
\begin{tabular}{c c | c c}
\hline
\hspace{0.5pc} Resolution, $N_{\iin}$
\hspace{0.5pc} & GMRES iterations 
\hspace{0.5pc} &  Number of bodies
\hspace{0.5pc} & GMRES iterations  \\
\hspace{0.0pc} (10 bodies) &
\hspace{0.5pc} & ($N_\iin = 128$ fixed) & \\
\hline
32	& 228	& 10 	& 223	\\
64    	& 226	& 15 	& 292	\\
128	& 223	& 20 	& 363	\\
256	& 223	& 30 	& 496	\\
512	& 223	& 40 	& 676	\\
1024	& 223	& 50 	& 991	\\
\hline
\end{tabular}
\end{center}
\end{table}

Lastly, we want to briefly examine the number of GMRES iterations required by our simulations since this is typically the most expensive step. Due to our {\em second-kind} boundary integral formulation, we expect the number of GMRES iterations to be independent of the spatial resolution, $N_{\iin}$. The iteration count, however, may depend on the number of bodies, as well as other factors such as the shapes of the bodies and their configuration. To examine a few of these factors, we show in Table~\ref{itertab} the number of GMRES iterations required to reach a specified tolerance of $10^{-8}$ in several different simulations. The left side of the table shows the results for a 10-body simulation, in which we keep the geometry fixed and vary the spatial resolution $N_{\iin}$. As seen in the table, the iteration count remains essentially fixed as the resolution increases. In the right side of the table, we keep the spatial resolution per body fixed and vary the number of simulated bodies from 10 to 50. The bodies are all circular with randomly assigned positions and radii (overlapping bodies are avoided). The table shows that the GMRES count increases roughly linearly with the number of bodies. 

\subsection{Multibody results}

\begin{figure}
\begin{center}
\includegraphics[width = 0.80 \textwidth]{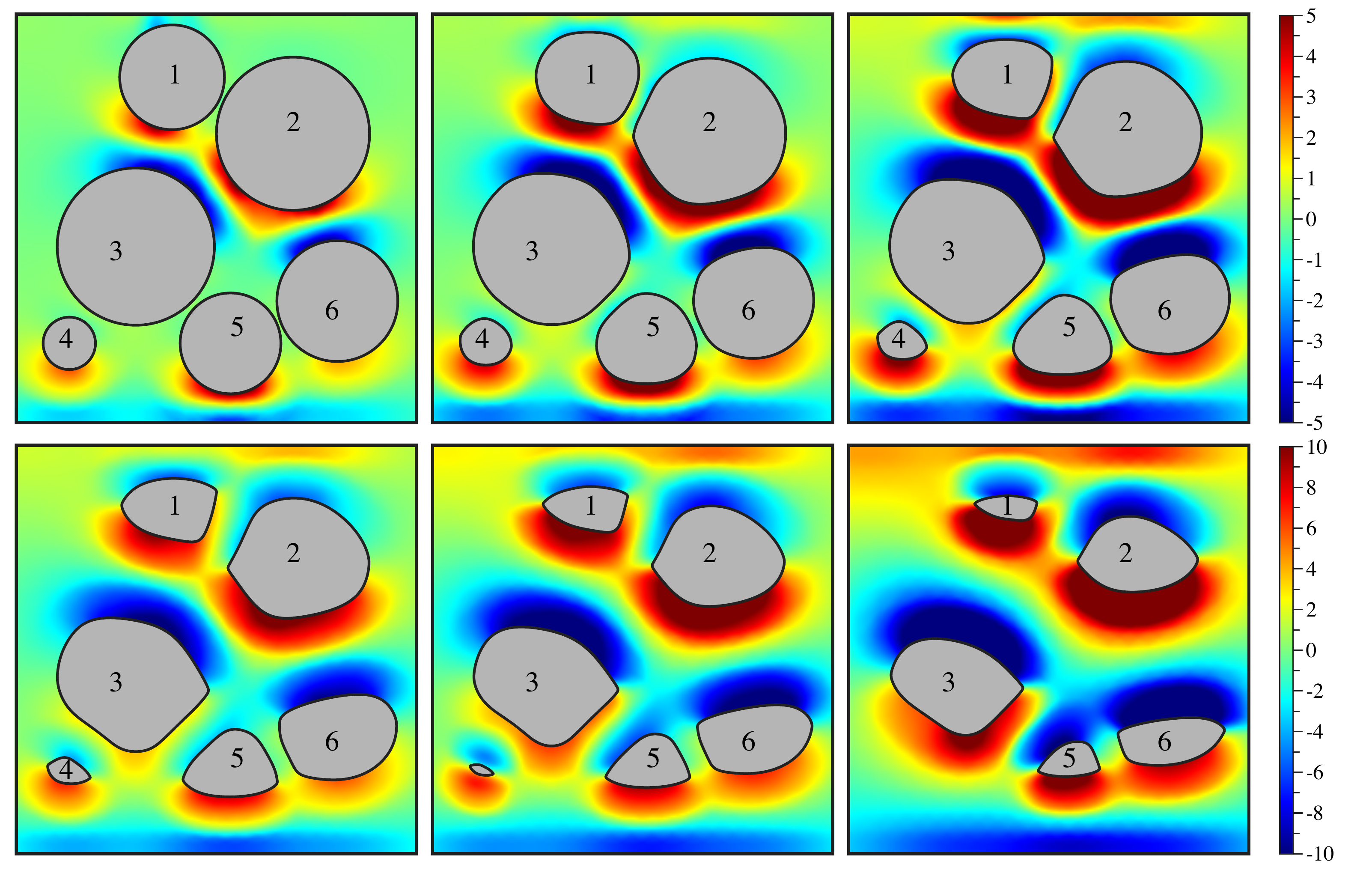}
\caption{\label{fig:06bodies} The simulated erosion of 6 bodies with the fixed-pressure-drop condition. Erosion causes flat faces to develop wherever bodies were initially very close to one another. This effect leads to the appearance of straight channels between bodies, and causes the bodies themselves to appear distinctly polygonal. The vorticity field (color-bar at right) shows that the flow is concentrated primarily in a single channel that runs horizontally across the domain. }
\end{center}
\end{figure}

With the validation complete, we now discuss the emergent physics of multibody erosion revealed by our simulations. We begin with a relatively simple case of 6 eroding bodies, as this will allow us to identify key principles at work in more complex cases (e.g.~the 50-body simulation in Fig.~\ref{fig:50bodies}). Figure~\ref{fig:06bodies} shows the simulated erosion of 6 bodies, each initial circles of various sizes and positions. Over time, erosion not only diminishes the size of the bodies, but also alters their shapes considerably. As in the single-body simulations, corner-like features develop, but not in locations that are easily predictable. Rather, when and where these corners form depends on how neighboring bodies interact with one another as mediated by the Stokes flow.

The initial configuration in  Fig.~\ref{fig:06bodies} was selected to have several points at which the bodies nearly contact one another, for example between bodies 1 \& 2, 3 \& 5, and 5 \& 6. Intriguingly, erosion appears to remove material preferentially from these sights, as the narrow spacing between neighboring bodies expands rather quickly. For example, the spacing between bodies 1 \& 2 is initially quite small, but by the second frame, is comparable to the spacing between bodies 3 \& 4. Moreover, the bodies tend to develop flat faces at the sites of near contact, leading to the appearance of straight channels between bodies. Straight channels can be seen, for example, between bodies 1 \& 2, 3 \& 5, and 5 \& 6. Once formed, these channels persist until the bodies have nearly vanished. Notice that by the final frame, body 4 has vanished completely, and the simulation handles this event without difficulty.

In Fig.~\ref{fig:06bodies}, we also show the vorticity field of the surrounding flow. This field helps not only to visualize regions of high shear, but also to identify certain preferred channels in which the flow is most intense.  As seen in the figure, the flow appears to be concentrated primarily within a single channel that runs horizontally below bodies 1 \& 2 and above bodies 3, 5, \& 6. This main channel widens the fastest, and secondary channels (e.g.~between bodies 1 \& 2 or 3 \& 5 or 5 \& 6) follow. Notice that the overall vorticity magnitude increases in time as the flow rate $\umax$ increases due to keeping the pressure drop fixed.

Looking back at the 50-body simulation in Fig.~\ref{fig:50bodies}, many of these same features can be observed. Bodies tend to develop corner-like features connected by flat faces, which causes them to appear polygonal. Nearly straight channels form between the bodies, and a few of these appear to be preferred channels which transmit the most flow and therefore widen the fastest.

\begin{figure}
\begin{center}
\includegraphics[width = 0.90 \textwidth]{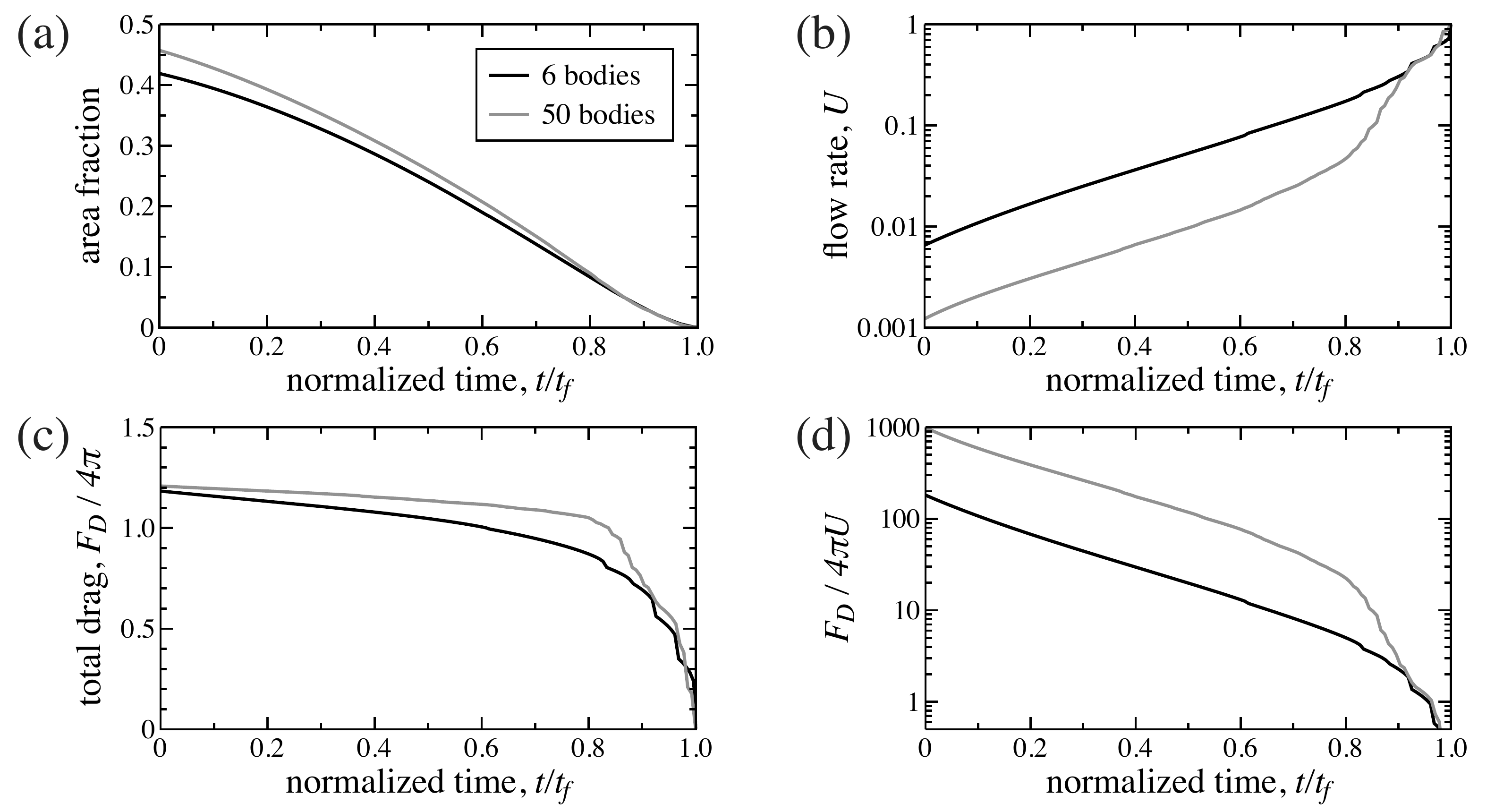}
\caption{\label{fig:mbodyplots} 
Multibody simulation results: black curves correspond to the 6-body simulation and gray to the 50-body simulation. (a) The area fraction of the eroding bodies versus time. (b) The incoming flow rate $\umax$ versus time. The fixed-pressure-drop condition causes $\umax$ to increase as bodies erode. (c) The total drag exerted by all bodies. Because $\umax$ increases as the bodies erode, the drag remains nearly constant throughout much of the simulation. In the very last moments, the drag drops suddenly to zero. (d) Normalizing the drag by $\umax$ provides gives an estimate of how much resistance the porous matrix provides. The resistance decreases rapidly as the bodies erode.
}
\end{center}
\end{figure}

To further quantify these observations, Fig.~\ref{fig:mbodyplots} shows several measurements taken from the multibody simulations. In all cases, the black curves correspond to the above 6-body simulation and gray curves to the 50-body simulation. First, Fig.~\ref{fig:mbodyplots}(a) shows the fraction of area occupied by the solid bodies as it varies in time. In the 6-body simulation, the initial configuration comprises about 40\% of the domain area, and this percentage steadily decreases as the bodies erode. The 50-body simulation shows a similar trend. 

Next, Figure~\ref{fig:mbodyplots}(b) shows how the flow rate, $\umax$, changes in time as a result of fixing the pressure drop across the cell. Notice that, in this figure, the vertical axis is logarithmic. Both simulations show the flow rate to increase rapidly, even during the early stages of erosion when the area fraction has not changed that significantly. Due to the logarithmic axis, the nearly linear trend seen in the figure corresponds to exponential growth of $\umax$. The change in area fraction alone cannot account for this rapid growth, suggesting the reshaping process to be vital. In particular, the formation of channels observed in Fig.~\ref{fig:06bodies} may account for these higher flow rates.

Naturally, the flow transmitted by the medium is linked to the amount of drag the bodies exert. We therefore show in Fig.~\ref{fig:mbodyplots}(c) the total drag exerted by all bodies as it varies in time. In both cases, the drag remains nearly constant throughout much of the simulation and then drops to zero in the final moments. Although all bodies are diminishing in size, the associated increase in $\umax$ leads to the nearly constant drag that is observed here.

To remove the effect of increasing flow rate, we show in Fig.~\ref{fig:mbodyplots}(d) the drag normalized by $\umax$. Since both the pressure and viscous drag depend linearly on $\umax$, this normalized quantity, $F_D / (4 \pi \umax)$, is {\em independent} of $\umax$ and therefore isolates how drag depends on the configuration of bodies. In this way, we can interpret $F_D / (4 \pi \umax)$ as the {\em resistance} associated with the porous network. As seen in the figure, the medium's resistance decreases rapidly as its constituents erode. Since the vertical axis is logarithmic, the rate is roughly exponential at early times. As before, the reduction in area-fraction is not sufficient to explain this rapid decline, leaving the channelization process as the most likely mechanism.

\section{Conclusions\label{s:conclusions}}

We have combined two accurate and stable numerical methods to simulate two-dimensional viscous erosion in a porous medium. First, a BIE formulation of the incompressible Stokes equations discretizes only the boundaries and resolves non-negligible interactions between bodies. This Stokes solver enables computation of the surface shear stress, which dictates the erosion rate of solid material. To advance solid boundaries forward in time, we use the {\thL} framework, which prevents tangling and distortion of the mesh that would otherwise lead to numerical instability. The resulting method is accurate (spectral in space, second-order in time), numerically stable, and accelerated via the fast multipole method.

Analytical results for the case of a single eroding body provide a stringent benchmark for our simulation. In particular, by appealing to the principle of long-time uniform stress, we predict a front/rear corner of angle $102^{\circ}$, which matches closely the value measured in the simulations. Scaling laws for the asymptotic rate of area and drag reduction also match the simulation.

Simulations of multiple eroding bodies reveals even more complex behavior. As in the single body case, we witness the formation of sharp corners, but not in locations that are easily predictable. Furthermore, bodies tend to form flat faces in places where there was initially near contact, leading to the appearance of straight channels between bodies. The overall resistance of the porous network declines rapidly with erosion, which we attribute to the channelization process.

In work currently underway, we are using this method to more thoroughly examine how bulk and statistical properties of the porous medium change with erosion. We are particularly interested in the development of anisotropic permeability of the medium. Preliminary simulations show promise. We furthermore would like to investigate how the statistical distribution of body size, shape, and alignment changes over time. We hope to report on these topics soon.

As a longer term goal, we would like to include the effect of sedimentation to more faithfully capture the physics of a porous medium.  Sedimentation requires overcoming a few challenges. First, the bodies must be mobile, as has been accomplished in previous work for dense concentrations of deformable and rigid bodies~\cite{qua-bir2014a, cor-gre-rac-vee2017, cor-vee2017, kli-tor2016b, rac-gre2016, qua-bir2016}. Second, contacts and overlap between bodies must be prevented. This second requirement has been recognized as quite a challenge, with many novel methods developed to treat it~\cite{lu-rah-zor2017, kab-qua-bir2018, san-mo1994}. Incorporation of one or more of these schemes, will enable simulations of a porous medium undergoing simultaneous erosion and sedimentation.

\paragraph{\bf Acknowledgments} We would like to thank Manas Rachh for
supplying the FMM for the Stokes double-layer potential. BQ and NM were
supported by Florida State University startup funds and Simons
Foundation Mathematics and Physical Sciences-Collaboration Grants for
Mathematicians.

\bibliographystyle{plainnat} 
\biboptions{sort&compress}


\end{document}

%% file: preambles.tex
\usepackage{amsfonts}
\usepackage[fleqn,reqno]{amsmath}
\usepackage{amssymb}
\usepackage[titletoc]{appendix}
\usepackage{array}
\usepackage{enumitem}
\usepackage{filecontents}
\usepackage[top=1.2in,bottom=1.2in,left=1in, right=1in]{geometry}
\usepackage{graphics}
\usepackage{lineno}
\usepackage{pgfplots}
\usepackage{tikz}
\usepackage{todonotes}
\usetikzlibrary{arrows}
\usepackage{comment}

\usepackage[pagebackref=false,bookmarks=false]{hyperref} 

\hypersetup{
  bookmarksnumbered=true,
  bookmarksopen=false,
  hypertexnames=false,      
  breaklinks=true,          
  unicode=false,
  pdffitwindow=true,        
  pdfnewwindow=true,        
  colorlinks=true,         
  linkcolor=dblue,
  anchorcolor=red,
  citecolor=dorange,
  filecolor=magenta,
  urlcolor=dblue,
  pdfstartview = FitH,
  pdfkeywords = {},
  pdfcreator = {LaTeX with hyperref package}
}

\newcommand{\bd}{{\partial}}
\newcommand{\bigO}{{\mathcal{O}}}
\newcommand{\cc}{{\mathbf{c}}}
\newcommand{\DD}{{\mathcal{D}}}
\newcommand{\eeta}{{\boldsymbol\eta}}
\newcommand{\ff}{{\mathbf{f}}}
\newcommand{\grad}{{\nabla}}
\newcommand{\II}{{\mathbf{I}}}
\newcommand{\iin}{\mathrm{in}}
\newcommand{\llambda}{{\boldsymbol\lambda}}
\newcommand{\nn}{{\mathbf{n}}}
\newcommand{\NN}{{\mathcal{N}}}
\newcommand{\out}{\mathrm{out}}
\newcommand{\rr}{{\mathbf{r}}}

\renewcommand{\ss}{{\mathbf{s}}}
\newcommand{\ssigma}{{\boldsymbol\sigma}}

\newcommand{\uu}{{\mathbf{u}}}
\newcommand{\UU}{{\mathbf{U}}}

\newcommand{\xx}{{\mathbf{x}}}

\newcommand{\yy}{{\mathbf{y}}}

\def\gap{\hspace*{.2in}}

\newcommand{\pderiv}[2]{\frac{\partial #1}{\partial #2}}
\newcommand{\tderiv}[2]{\frac{d #1}{d #2}}
\newcommand{\ppd}[2]{\frac{\partial^2 #1}{{\partial #2}^2}}
\newcommand{\pdi}[2]{\partial {#1}/\partial {#2}}

\newcommand{\abs}[1]{\lvert #1 \rvert}
\newcommand{\mean}[1]{\left< #1 \right>}
\newcommand{\thL}{$\theta$--$L$}
\newcommand{\eps}{\varepsilon}
\newcommand{\Vn}{V_\nn}
\newcommand{\Vs}{V_\ss}
\newcommand{\atau}{\abs{\tau}}
\newcommand{\thalpha}{\pderiv{\theta}{\alpha}}
\newcommand{\elfun}{\zeta}
\newcommand{\thhat}{\hat{\theta}}
\newcommand{\Dt}{\Delta t}
\newcommand{\NLterm}{\mathcal{N}}
\newcommand{\Mterm}{\mathcal{M}}
\newcommand{\FourierSum}{ \sum_{k = -N_\iin /2}^{N_\iin /2-1} }
\newcommand{\atausig}{\atau^{(\sigma)}}
\newcommand{\Vnsig}{\Vn^{(\sigma)}}
\newcommand{\Vssig}{\Vs^{(\sigma)}}
\newcommand{\umax}{U}

\newcommand{\FFD}{\mathbf{F}_{D}}
\newcommand{\oangle}{\beta}
\newcommand{\thb}{\theta_0}

%% file: schematic.tikz
\begin{tikzpicture}[scale=1.5] 

\begin{axis}[ 
axis equal image, 
scale only axis, 
xmin=-3.04, 
xmax=3.6, 
ymin=-1.1, 
ymax=1.1, 
hide axis, 
] 

\addplot [color=black,dashed,line width=1] coordinates{ 
  (-1,-1)
  (-1,+1)
}; 

\addplot [color=black,dashed,line width=1] coordinates{ 
  (1,-1)
  (1,+1)
}; 

\addplot [color=black,solid,line width=2] coordinates{ 
(3.0000e+00,0.0000e+00)
(3.0000e+00,2.4549e-02)
(3.0000e+00,4.9127e-02)
(3.0000e+00,7.3764e-02)
(3.0000e+00,9.8491e-02)
(3.0000e+00,1.2334e-01)
(3.0000e+00,1.4834e-01)
(3.0000e+00,1.7352e-01)
(3.0000e+00,1.9891e-01)
(3.0000e+00,2.2456e-01)
(3.0000e+00,2.5049e-01)
(3.0000e+00,2.7674e-01)
(3.0000e+00,3.0334e-01)
(2.9999e+00,3.3035e-01)
(2.9999e+00,3.5779e-01)
(2.9998e+00,3.8572e-01)
(2.9997e+00,4.1417e-01)
(2.9994e+00,4.4319e-01)
(2.9991e+00,4.7282e-01)
(2.9985e+00,5.0310e-01)
(2.9975e+00,5.3407e-01)
(2.9960e+00,5.6575e-01)
(2.9938e+00,5.9814e-01)
(2.9904e+00,6.3123e-01)
(2.9854e+00,6.6493e-01)
(2.9780e+00,6.9912e-01)
(2.9673e+00,7.3358e-01)
(2.9520e+00,7.6793e-01)
(2.9306e+00,8.0171e-01)
(2.9014e+00,8.3425e-01)
(2.8625e+00,8.6481e-01)
(2.8126e+00,8.9262e-01)
(2.7510e+00,9.1700e-01)
(2.6779e+00,9.3755e-01)
(2.5944e+00,9.5418e-01)
(2.5028e+00,9.6713e-01)
(2.4051e+00,9.7688e-01)
(2.3038e+00,9.8402e-01)
(2.2007e+00,9.8911e-01)
(2.0974e+00,9.9268e-01)
(1.9948e+00,9.9514e-01)
(1.8937e+00,9.9681e-01)
(1.7944e+00,9.9794e-01)
(1.6972e+00,9.9868e-01)
(1.6022e+00,9.9917e-01)
(1.5093e+00,9.9949e-01)
(1.4185e+00,9.9969e-01)
(1.3296e+00,9.9981e-01)
(1.2425e+00,9.9989e-01)
(1.1572e+00,9.9994e-01)
(1.0734e+00,9.9997e-01)
(9.9105e-01,9.9998e-01)
(9.1003e-01,9.9999e-01)
(8.3021e-01,1.0000e+00)
(7.5146e-01,1.0000e+00)
(6.7367e-01,1.0000e+00)
(5.9674e-01,1.0000e+00)
(5.2055e-01,1.0000e+00)
(4.4501e-01,1.0000e+00)
(3.7001e-01,1.0000e+00)
(2.9547e-01,1.0000e+00)
(2.2129e-01,1.0000e+00)
(1.4738e-01,1.0000e+00)
(7.3646e-02,1.0000e+00)
(1.8370e-16,1.0000e+00)
(-7.3646e-02,1.0000e+00)
(-1.4738e-01,1.0000e+00)
(-2.2129e-01,1.0000e+00)
(-2.9547e-01,1.0000e+00)
(-3.7001e-01,1.0000e+00)
(-4.4501e-01,1.0000e+00)
(-5.2055e-01,1.0000e+00)
(-5.9674e-01,1.0000e+00)
(-6.7367e-01,1.0000e+00)
(-7.5146e-01,1.0000e+00)
(-8.3021e-01,1.0000e+00)
(-9.1003e-01,9.9999e-01)
(-9.9105e-01,9.9998e-01)
(-1.0734e+00,9.9997e-01)
(-1.1572e+00,9.9994e-01)
(-1.2425e+00,9.9989e-01)
(-1.3296e+00,9.9981e-01)
(-1.4185e+00,9.9969e-01)
(-1.5093e+00,9.9949e-01)
(-1.6022e+00,9.9917e-01)
(-1.6972e+00,9.9868e-01)
(-1.7944e+00,9.9794e-01)
(-1.8937e+00,9.9681e-01)
(-1.9948e+00,9.9514e-01)
(-2.0974e+00,9.9268e-01)
(-2.2007e+00,9.8911e-01)
(-2.3038e+00,9.8402e-01)
(-2.4051e+00,9.7688e-01)
(-2.5028e+00,9.6713e-01)
(-2.5944e+00,9.5418e-01)
(-2.6779e+00,9.3755e-01)
(-2.7510e+00,9.1700e-01)
(-2.8126e+00,8.9262e-01)
(-2.8625e+00,8.6481e-01)
(-2.9014e+00,8.3425e-01)
(-2.9306e+00,8.0171e-01)
(-2.9520e+00,7.6793e-01)
(-2.9673e+00,7.3358e-01)
(-2.9780e+00,6.9912e-01)
(-2.9854e+00,6.6493e-01)
(-2.9904e+00,6.3123e-01)
(-2.9938e+00,5.9814e-01)
(-2.9960e+00,5.6575e-01)
(-2.9975e+00,5.3407e-01)
(-2.9985e+00,5.0310e-01)
(-2.9991e+00,4.7282e-01)
(-2.9994e+00,4.4319e-01)
(-2.9997e+00,4.1417e-01)
(-2.9998e+00,3.8572e-01)
(-2.9999e+00,3.5779e-01)
(-2.9999e+00,3.3035e-01)
(-3.0000e+00,3.0334e-01)
(-3.0000e+00,2.7674e-01)
(-3.0000e+00,2.5049e-01)
(-3.0000e+00,2.2456e-01)
(-3.0000e+00,1.9891e-01)
(-3.0000e+00,1.7352e-01)
(-3.0000e+00,1.4834e-01)
(-3.0000e+00,1.2334e-01)
(-3.0000e+00,9.8491e-02)
(-3.0000e+00,7.3764e-02)
(-3.0000e+00,4.9127e-02)
(-3.0000e+00,2.4549e-02)
(-3.0000e+00,1.2246e-16)
(-3.0000e+00,-2.4549e-02)
(-3.0000e+00,-4.9127e-02)
(-3.0000e+00,-7.3764e-02)
(-3.0000e+00,-9.8491e-02)
(-3.0000e+00,-1.2334e-01)
(-3.0000e+00,-1.4834e-01)
(-3.0000e+00,-1.7352e-01)
(-3.0000e+00,-1.9891e-01)
(-3.0000e+00,-2.2456e-01)
(-3.0000e+00,-2.5049e-01)
(-3.0000e+00,-2.7674e-01)
(-3.0000e+00,-3.0334e-01)
(-2.9999e+00,-3.3035e-01)
(-2.9999e+00,-3.5779e-01)
(-2.9998e+00,-3.8572e-01)
(-2.9997e+00,-4.1417e-01)
(-2.9994e+00,-4.4319e-01)
(-2.9991e+00,-4.7282e-01)
(-2.9985e+00,-5.0310e-01)
(-2.9975e+00,-5.3407e-01)
(-2.9960e+00,-5.6575e-01)
(-2.9938e+00,-5.9814e-01)
(-2.9904e+00,-6.3123e-01)
(-2.9854e+00,-6.6493e-01)
(-2.9780e+00,-6.9912e-01)
(-2.9673e+00,-7.3358e-01)
(-2.9520e+00,-7.6793e-01)
(-2.9306e+00,-8.0171e-01)
(-2.9014e+00,-8.3425e-01)
(-2.8625e+00,-8.6481e-01)
(-2.8126e+00,-8.9262e-01)
(-2.7510e+00,-9.1700e-01)
(-2.6779e+00,-9.3755e-01)
(-2.5944e+00,-9.5418e-01)
(-2.5028e+00,-9.6713e-01)
(-2.4051e+00,-9.7688e-01)
(-2.3038e+00,-9.8402e-01)
(-2.2007e+00,-9.8911e-01)
(-2.0974e+00,-9.9268e-01)
(-1.9948e+00,-9.9514e-01)
(-1.8937e+00,-9.9681e-01)
(-1.7944e+00,-9.9794e-01)
(-1.6972e+00,-9.9868e-01)
(-1.6022e+00,-9.9917e-01)
(-1.5093e+00,-9.9949e-01)
(-1.4185e+00,-9.9969e-01)
(-1.3296e+00,-9.9981e-01)
(-1.2425e+00,-9.9989e-01)
(-1.1572e+00,-9.9994e-01)
(-1.0734e+00,-9.9997e-01)
(-9.9105e-01,-9.9998e-01)
(-9.1003e-01,-9.9999e-01)
(-8.3021e-01,-1.0000e+00)
(-7.5146e-01,-1.0000e+00)
(-6.7367e-01,-1.0000e+00)
(-5.9674e-01,-1.0000e+00)
(-5.2055e-01,-1.0000e+00)
(-4.4501e-01,-1.0000e+00)
(-3.7001e-01,-1.0000e+00)
(-2.9547e-01,-1.0000e+00)
(-2.2129e-01,-1.0000e+00)
(-1.4738e-01,-1.0000e+00)
(-7.3646e-02,-1.0000e+00)
(-5.5109e-16,-1.0000e+00)
(7.3646e-02,-1.0000e+00)
(1.4738e-01,-1.0000e+00)
(2.2129e-01,-1.0000e+00)
(2.9547e-01,-1.0000e+00)
(3.7001e-01,-1.0000e+00)
(4.4501e-01,-1.0000e+00)
(5.2055e-01,-1.0000e+00)
(5.9674e-01,-1.0000e+00)
(6.7367e-01,-1.0000e+00)
(7.5146e-01,-1.0000e+00)
(8.3021e-01,-1.0000e+00)
(9.1003e-01,-9.9999e-01)
(9.9105e-01,-9.9998e-01)
(1.0734e+00,-9.9997e-01)
(1.1572e+00,-9.9994e-01)
(1.2425e+00,-9.9989e-01)
(1.3296e+00,-9.9981e-01)
(1.4185e+00,-9.9969e-01)
(1.5093e+00,-9.9949e-01)
(1.6022e+00,-9.9917e-01)
(1.6972e+00,-9.9868e-01)
(1.7944e+00,-9.9794e-01)
(1.8937e+00,-9.9681e-01)
(1.9948e+00,-9.9514e-01)
(2.0974e+00,-9.9268e-01)
(2.2007e+00,-9.8911e-01)
(2.3038e+00,-9.8402e-01)
(2.4051e+00,-9.7688e-01)
(2.5028e+00,-9.6713e-01)
(2.5944e+00,-9.5418e-01)
(2.6779e+00,-9.3755e-01)
(2.7510e+00,-9.1700e-01)
(2.8126e+00,-8.9262e-01)
(2.8625e+00,-8.6481e-01)
(2.9014e+00,-8.3425e-01)
(2.9306e+00,-8.0171e-01)
(2.9520e+00,-7.6793e-01)
(2.9673e+00,-7.3358e-01)
(2.9780e+00,-6.9912e-01)
(2.9854e+00,-6.6493e-01)
(2.9904e+00,-6.3123e-01)
(2.9938e+00,-5.9814e-01)
(2.9960e+00,-5.6575e-01)
(2.9975e+00,-5.3407e-01)
(2.9985e+00,-5.0310e-01)
(2.9991e+00,-4.7282e-01)
(2.9994e+00,-4.4319e-01)
(2.9997e+00,-4.1417e-01)
(2.9998e+00,-3.8572e-01)
(2.9999e+00,-3.5779e-01)
(2.9999e+00,-3.3035e-01)
(3.0000e+00,-3.0334e-01)
(3.0000e+00,-2.7674e-01)
(3.0000e+00,-2.5049e-01)
(3.0000e+00,-2.2456e-01)
(3.0000e+00,-1.9891e-01)
(3.0000e+00,-1.7352e-01)
(3.0000e+00,-1.4834e-01)
(3.0000e+00,-1.2334e-01)
(3.0000e+00,-9.8491e-02)
(3.0000e+00,-7.3764e-02)
(3.0000e+00,-4.9127e-02)
(3.0000e+00,-2.4549e-02)
(3.0000e+00,0.0000e+00)
}; 

\addplot [color=black,solid,fill] coordinates{ 
(5.0000e-01,0.0000e+00)
(4.9759e-01,4.9009e-02)
(4.9039e-01,9.7545e-02)
(4.7847e-01,1.4514e-01)
(4.6194e-01,1.9134e-01)
(4.4096e-01,2.3570e-01)
(4.1573e-01,2.7779e-01)
(3.8651e-01,3.1720e-01)
(3.5355e-01,3.5355e-01)
(3.1720e-01,3.8651e-01)
(2.7779e-01,4.1573e-01)
(2.3570e-01,4.4096e-01)
(1.9134e-01,4.6194e-01)
(1.4514e-01,4.7847e-01)
(9.7545e-02,4.9039e-01)
(4.9009e-02,4.9759e-01)
(3.0616e-17,5.0000e-01)
(-4.9009e-02,4.9759e-01)
(-9.7545e-02,4.9039e-01)
(-1.4514e-01,4.7847e-01)
(-1.9134e-01,4.6194e-01)
(-2.3570e-01,4.4096e-01)
(-2.7779e-01,4.1573e-01)
(-3.1720e-01,3.8651e-01)
(-3.5355e-01,3.5355e-01)
(-3.8651e-01,3.1720e-01)
(-4.1573e-01,2.7779e-01)
(-4.4096e-01,2.3570e-01)
(-4.6194e-01,1.9134e-01)
(-4.7847e-01,1.4514e-01)
(-4.9039e-01,9.7545e-02)
(-4.9759e-01,4.9009e-02)
(-5.0000e-01,6.1232e-17)
(-4.9759e-01,-4.9009e-02)
(-4.9039e-01,-9.7545e-02)
(-4.7847e-01,-1.4514e-01)
(-4.6194e-01,-1.9134e-01)
(-4.4096e-01,-2.3570e-01)
(-4.1573e-01,-2.7779e-01)
(-3.8651e-01,-3.1720e-01)
(-3.5355e-01,-3.5355e-01)
(-3.1720e-01,-3.8651e-01)
(-2.7779e-01,-4.1573e-01)
(-2.3570e-01,-4.4096e-01)
(-1.9134e-01,-4.6194e-01)
(-1.4514e-01,-4.7847e-01)
(-9.7545e-02,-4.9039e-01)
(-4.9009e-02,-4.9759e-01)
(-9.1849e-17,-5.0000e-01)
(4.9009e-02,-4.9759e-01)
(9.7545e-02,-4.9039e-01)
(1.4514e-01,-4.7847e-01)
(1.9134e-01,-4.6194e-01)
(2.3570e-01,-4.4096e-01)
(2.7779e-01,-4.1573e-01)
(3.1720e-01,-3.8651e-01)
(3.5355e-01,-3.5355e-01)
(3.8651e-01,-3.1720e-01)
(4.1573e-01,-2.7779e-01)
(4.4096e-01,-2.3570e-01)
(4.6194e-01,-1.9134e-01)
(4.7847e-01,-1.4514e-01)
(4.9039e-01,-9.7545e-02)
(4.9759e-01,-4.9009e-02)
(5.0000e-01,-1.2246e-16)
};

\addplot [color=black,solid,fill] coordinates{ 
(8.0000e-01,5.0000e-01)
(7.9952e-01,5.0980e-01)
(7.9808e-01,5.1951e-01)
(7.9569e-01,5.2903e-01)
(7.9239e-01,5.3827e-01)
(7.8819e-01,5.4714e-01)
(7.8315e-01,5.5556e-01)
(7.7730e-01,5.6344e-01)
(7.7071e-01,5.7071e-01)
(7.6344e-01,5.7730e-01)
(7.5556e-01,5.8315e-01)
(7.4714e-01,5.8819e-01)
(7.3827e-01,5.9239e-01)
(7.2903e-01,5.9569e-01)
(7.1951e-01,5.9808e-01)
(7.0980e-01,5.9952e-01)
(7.0000e-01,6.0000e-01)
(6.9020e-01,5.9952e-01)
(6.8049e-01,5.9808e-01)
(6.7097e-01,5.9569e-01)
(6.6173e-01,5.9239e-01)
(6.5286e-01,5.8819e-01)
(6.4444e-01,5.8315e-01)
(6.3656e-01,5.7730e-01)
(6.2929e-01,5.7071e-01)
(6.2270e-01,5.6344e-01)
(6.1685e-01,5.5556e-01)
(6.1181e-01,5.4714e-01)
(6.0761e-01,5.3827e-01)
(6.0431e-01,5.2903e-01)
(6.0192e-01,5.1951e-01)
(6.0048e-01,5.0980e-01)
(6.0000e-01,5.0000e-01)
(6.0048e-01,4.9020e-01)
(6.0192e-01,4.8049e-01)
(6.0431e-01,4.7097e-01)
(6.0761e-01,4.6173e-01)
(6.1181e-01,4.5286e-01)
(6.1685e-01,4.4444e-01)
(6.2270e-01,4.3656e-01)
(6.2929e-01,4.2929e-01)
(6.3656e-01,4.2270e-01)
(6.4444e-01,4.1685e-01)
(6.5286e-01,4.1181e-01)
(6.6173e-01,4.0761e-01)
(6.7097e-01,4.0431e-01)
(6.8049e-01,4.0192e-01)
(6.9020e-01,4.0048e-01)
(7.0000e-01,4.0000e-01)
(7.0980e-01,4.0048e-01)
(7.1951e-01,4.0192e-01)
(7.2903e-01,4.0431e-01)
(7.3827e-01,4.0761e-01)
(7.4714e-01,4.1181e-01)
(7.5556e-01,4.1685e-01)
(7.6344e-01,4.2270e-01)
(7.7071e-01,4.2929e-01)
(7.7730e-01,4.3656e-01)
(7.8315e-01,4.4444e-01)
(7.8819e-01,4.5286e-01)
(7.9239e-01,4.6173e-01)
(7.9569e-01,4.7097e-01)
(7.9808e-01,4.8049e-01)
(7.9952e-01,4.9020e-01)
(8.0000e-01,5.0000e-01)
};

\addplot [color=black,solid,fill] coordinates{ 
(-6.0000e-01,6.0000e-01)
(-6.0096e-01,6.1960e-01)
(-6.0384e-01,6.3902e-01)
(-6.0861e-01,6.5806e-01)
(-6.1522e-01,6.7654e-01)
(-6.2362e-01,6.9428e-01)
(-6.3371e-01,7.1111e-01)
(-6.4540e-01,7.2688e-01)
(-6.5858e-01,7.4142e-01)
(-6.7312e-01,7.5460e-01)
(-6.8889e-01,7.6629e-01)
(-7.0572e-01,7.7638e-01)
(-7.2346e-01,7.8478e-01)
(-7.4194e-01,7.9139e-01)
(-7.6098e-01,7.9616e-01)
(-7.8040e-01,7.9904e-01)
(-8.0000e-01,8.0000e-01)
(-8.1960e-01,7.9904e-01)
(-8.3902e-01,7.9616e-01)
(-8.5806e-01,7.9139e-01)
(-8.7654e-01,7.8478e-01)
(-8.9428e-01,7.7638e-01)
(-9.1111e-01,7.6629e-01)
(-9.2688e-01,7.5460e-01)
(-9.4142e-01,7.4142e-01)
(-9.5460e-01,7.2688e-01)
(-9.6629e-01,7.1111e-01)
(-9.7638e-01,6.9428e-01)
(-9.8478e-01,6.7654e-01)
(-9.9139e-01,6.5806e-01)
(-9.9616e-01,6.3902e-01)
(-9.9904e-01,6.1960e-01)
(-1.0000e+00,6.0000e-01)
(-9.9904e-01,5.8040e-01)
(-9.9616e-01,5.6098e-01)
(-9.9139e-01,5.4194e-01)
(-9.8478e-01,5.2346e-01)
(-9.7638e-01,5.0572e-01)
(-9.6629e-01,4.8889e-01)
(-9.5460e-01,4.7312e-01)
(-9.4142e-01,4.5858e-01)
(-9.2688e-01,4.4540e-01)
(-9.1111e-01,4.3371e-01)
(-8.9428e-01,4.2362e-01)
(-8.7654e-01,4.1522e-01)
(-8.5806e-01,4.0861e-01)
(-8.3902e-01,4.0384e-01)
(-8.1960e-01,4.0096e-01)
(-8.0000e-01,4.0000e-01)
(-7.8040e-01,4.0096e-01)
(-7.6098e-01,4.0384e-01)
(-7.4194e-01,4.0861e-01)
(-7.2346e-01,4.1522e-01)
(-7.0572e-01,4.2362e-01)
(-6.8889e-01,4.3371e-01)
(-6.7312e-01,4.4540e-01)
(-6.5858e-01,4.5858e-01)
(-6.4540e-01,4.7312e-01)
(-6.3371e-01,4.8889e-01)
(-6.2362e-01,5.0572e-01)
(-6.1522e-01,5.2346e-01)
(-6.0861e-01,5.4194e-01)
(-6.0384e-01,5.6098e-01)
(-6.0096e-01,5.8040e-01)
(-6.0000e-01,6.0000e-01)
};

\addplot [color=black,solid,fill] coordinates{ 
(9.0000e-01,-6.0000e-01)
(8.9904e-01,-5.8040e-01)
(8.9616e-01,-5.6098e-01)
(8.9139e-01,-5.4194e-01)
(8.8478e-01,-5.2346e-01)
(8.7638e-01,-5.0572e-01)
(8.6629e-01,-4.8889e-01)
(8.5460e-01,-4.7312e-01)
(8.4142e-01,-4.5858e-01)
(8.2688e-01,-4.4540e-01)
(8.1111e-01,-4.3371e-01)
(7.9428e-01,-4.2362e-01)
(7.7654e-01,-4.1522e-01)
(7.5806e-01,-4.0861e-01)
(7.3902e-01,-4.0384e-01)
(7.1960e-01,-4.0096e-01)
(7.0000e-01,-4.0000e-01)
(6.8040e-01,-4.0096e-01)
(6.6098e-01,-4.0384e-01)
(6.4194e-01,-4.0861e-01)
(6.2346e-01,-4.1522e-01)
(6.0572e-01,-4.2362e-01)
(5.8889e-01,-4.3371e-01)
(5.7312e-01,-4.4540e-01)
(5.5858e-01,-4.5858e-01)
(5.4540e-01,-4.7312e-01)
(5.3371e-01,-4.8889e-01)
(5.2362e-01,-5.0572e-01)
(5.1522e-01,-5.2346e-01)
(5.0861e-01,-5.4194e-01)
(5.0384e-01,-5.6098e-01)
(5.0096e-01,-5.8040e-01)
(5.0000e-01,-6.0000e-01)
(5.0096e-01,-6.1960e-01)
(5.0384e-01,-6.3902e-01)
(5.0861e-01,-6.5806e-01)
(5.1522e-01,-6.7654e-01)
(5.2362e-01,-6.9428e-01)
(5.3371e-01,-7.1111e-01)
(5.4540e-01,-7.2688e-01)
(5.5858e-01,-7.4142e-01)
(5.7312e-01,-7.5460e-01)
(5.8889e-01,-7.6629e-01)
(6.0572e-01,-7.7638e-01)
(6.2346e-01,-7.8478e-01)
(6.4194e-01,-7.9139e-01)
(6.6098e-01,-7.9616e-01)
(6.8040e-01,-7.9904e-01)
(7.0000e-01,-8.0000e-01)
(7.1960e-01,-7.9904e-01)
(7.3902e-01,-7.9616e-01)
(7.5806e-01,-7.9139e-01)
(7.7654e-01,-7.8478e-01)
(7.9428e-01,-7.7638e-01)
(8.1111e-01,-7.6629e-01)
(8.2688e-01,-7.5460e-01)
(8.4142e-01,-7.4142e-01)
(8.5460e-01,-7.2688e-01)
(8.6629e-01,-7.1111e-01)
(8.7638e-01,-6.9428e-01)
(8.8478e-01,-6.7654e-01)
(8.9139e-01,-6.5806e-01)
(8.9616e-01,-6.3902e-01)
(8.9904e-01,-6.1960e-01)
(9.0000e-01,-6.0000e-01)
};

\addplot [color=black,solid,fill] coordinates{ 
(-6.4000e-01,-5.0000e-01)
(-6.4029e-01,-4.9412e-01)
(-6.4115e-01,-4.8829e-01)
(-6.4258e-01,-4.8258e-01)
(-6.4457e-01,-4.7704e-01)
(-6.4708e-01,-4.7172e-01)
(-6.5011e-01,-4.6667e-01)
(-6.5362e-01,-4.6194e-01)
(-6.5757e-01,-4.5757e-01)
(-6.6194e-01,-4.5362e-01)
(-6.6667e-01,-4.5011e-01)
(-6.7172e-01,-4.4708e-01)
(-6.7704e-01,-4.4457e-01)
(-6.8258e-01,-4.4258e-01)
(-6.8829e-01,-4.4115e-01)
(-6.9412e-01,-4.4029e-01)
(-7.0000e-01,-4.4000e-01)
(-7.0588e-01,-4.4029e-01)
(-7.1171e-01,-4.4115e-01)
(-7.1742e-01,-4.4258e-01)
(-7.2296e-01,-4.4457e-01)
(-7.2828e-01,-4.4708e-01)
(-7.3333e-01,-4.5011e-01)
(-7.3806e-01,-4.5362e-01)
(-7.4243e-01,-4.5757e-01)
(-7.4638e-01,-4.6194e-01)
(-7.4989e-01,-4.6667e-01)
(-7.5292e-01,-4.7172e-01)
(-7.5543e-01,-4.7704e-01)
(-7.5742e-01,-4.8258e-01)
(-7.5885e-01,-4.8829e-01)
(-7.5971e-01,-4.9412e-01)
(-7.6000e-01,-5.0000e-01)
(-7.5971e-01,-5.0588e-01)
(-7.5885e-01,-5.1171e-01)
(-7.5742e-01,-5.1742e-01)
(-7.5543e-01,-5.2296e-01)
(-7.5292e-01,-5.2828e-01)
(-7.4989e-01,-5.3333e-01)
(-7.4638e-01,-5.3806e-01)
(-7.4243e-01,-5.4243e-01)
(-7.3806e-01,-5.4638e-01)
(-7.3333e-01,-5.4989e-01)
(-7.2828e-01,-5.5292e-01)
(-7.2296e-01,-5.5543e-01)
(-7.1742e-01,-5.5742e-01)
(-7.1171e-01,-5.5885e-01)
(-7.0588e-01,-5.5971e-01)
(-7.0000e-01,-5.6000e-01)
(-6.9412e-01,-5.5971e-01)
(-6.8829e-01,-5.5885e-01)
(-6.8258e-01,-5.5742e-01)
(-6.7704e-01,-5.5543e-01)
(-6.7172e-01,-5.5292e-01)
(-6.6667e-01,-5.4989e-01)
(-6.6194e-01,-5.4638e-01)
(-6.5757e-01,-5.4243e-01)
(-6.5362e-01,-5.3806e-01)
(-6.5011e-01,-5.3333e-01)
(-6.4708e-01,-5.2828e-01)
(-6.4457e-01,-5.2296e-01)
(-6.4258e-01,-5.1742e-01)
(-6.4115e-01,-5.1171e-01)
(-6.4029e-01,-5.0588e-01)
(-6.4000e-01,-5.0000e-01)
};


\node[font = \Huge,color=black] at (500,110) {$\Omega$};
\node[font = \normalsize,color=red] at (338,122) {$\nn$};
\node[font = \normalsize,color=red] at (342,158) {$\ss$};
\path[->,line width=1.2](240,90) edge (260,100);
\node at (225,90) {$\gamma_\ell$};
\path[->,line width=1.2,color=black](140,50) edge (120,14);
\node[font = \Large,color=black] at (150,70) {$\Gamma$};

\foreach \y in {-0.7,-0.5,...,0.7}
\addplot[color=black,line width = 1.0pt,solid,->]
plot coordinates{
  (-3,\y)
  (-3+0.6*(1-\y*\y),\y)
};

\foreach \y in {-0.7,-0.5,...,0.7}
\addplot[color=black,line width = 1.0pt,solid,->]
plot coordinates{
  (3,\y)
  (3+0.6*(1-\y*\y),\y)
};

\addplot[color=red,line width=0.5pt,solid,->]
plot coordinates{
  (0.3865,0.3172)
  (0.0966,0.0793)
};

\addplot[color=red,line width=0.5pt,solid,->]
plot coordinates{
  (0.3865,0.3172)
  (0.1486,0.6071)
};

\end{axis}

\end{tikzpicture}

%% file: erosion_methods.bbl
\begin{thebibliography}{72}
\providecommand{\natexlab}[1]{#1}
\providecommand{\url}[1]{\texttt{#1}}
\expandafter\ifx\csname urlstyle\endcsname\relax
  \providecommand{\doi}[1]{doi: #1}\else
  \providecommand{\doi}{doi: \begingroup \urlstyle{rm}\Url}\fi

\bibitem[Abrams et~al.(2009)Abrams, Lobkovsky, Petroff, Straub, McElroy,
  Mohrig, Kudrolli, and Rothman]{abrams2009growth}
Daniel~M Abrams, Alexander~E Lobkovsky, Alexander~P Petroff, Kyle~M Straub,
  Brandon McElroy, David~C Mohrig, Arshad Kudrolli, and Daniel~H Rothman.
\newblock Growth laws for channel networks incised by groundwater flow.
\newblock \emph{Nature Geoscience}, 2\penalty0 (3):\penalty0 193, 2009.

\bibitem[af~Klinteberg and Tornberg(2016)]{kli-tor2016b}
Ludvig af~Klinteberg and Anna-Karin Tornberg.
\newblock A fast integral equation method for solid particles in viscous flow
  using quadrature by expansion.
\newblock \emph{Journal of Computational Physics}, 326:\penalty0 420--445,
  2016.

\bibitem[Barnett et~al.(2015)Barnett, Wu, and Veerapaneni]{bar-wu-vee2015}
Alex Barnett, Bowei Wu, and Shravan Veerapaneni.
\newblock {Spectrally-Accurate Quadratures for Evaluation of Layer Potentials
  Close to the Boundary for the 2D Stokes and Laplace Equations}.
\newblock \emph{SIAM Journal on Scientific Computing}, 37\penalty0
  (4):\penalty0 B519--B542, 2015.

\bibitem[Barnett et~al.(2018)Barnett, Marple, Veerapaneni, and
  Zhao]{bar-mar-vee-zha2018}
Alex Barnett, Gary Marple, Shravan Veerapaneni, and Lin Zhao.
\newblock {A unified integral equation scheme for doubly-periodic Laplace and
  Stokes boundary value problems in two dimensions}.
\newblock \emph{Communications on Pure and Applied Mathematic}, 2018.
\newblock In press.

\bibitem[Barnett(2014)]{bar2014}
Alex~H. Barnett.
\newblock {Evaluation of layer potentials close to the boundary for Laplace and
  Helmholtz problems on analytic planar domains}.
\newblock \emph{SIAM Journal on Scientific Computing}, 36\penalty0
  (2):\penalty0 A427--A451, 2014.

\bibitem[Batchelor(1970)]{batchelor1970slender}
GK~Batchelor.
\newblock {Slender-body theory for particles of arbitrary cross-section in
  Stokes flow}.
\newblock \emph{Journal of Fluid Mechanics}, 44\penalty0 (3):\penalty0
  419--440, 1970.

\bibitem[Beale and Lai(2001)]{bea-lai2001}
J.T. Beale and M.-C. Lai.
\newblock {A Method for Computing Nearly Singular Integrals}.
\newblock \emph{SIAM Journal on Numerical Analysis}, 38\penalty0 (6):\penalty0
  1902--1925, 2001.

\bibitem[Bear(2013)]{bear2013dynamics}
Jacob Bear.
\newblock \emph{Dynamics of fluids in porous media}.
\newblock Courier Corporation, 2013.

\bibitem[Berhanu et~al.(2012)Berhanu, Petroff, Devauchelle, Kudrolli, and
  Rothman]{berhanu2012shape}
Michael Berhanu, Alexander Petroff, Olivier Devauchelle, Arshad Kudrolli, and
  Daniel~H Rothman.
\newblock Shape and dynamics of seepage erosion in a horizontal granular bed.
\newblock \emph{Physical Review E}, 86\penalty0 (4):\penalty0 041304, 2012.

\bibitem[Bremer(2012)]{bre2012}
James Bremer.
\newblock {On the Nystr\"om discretization of integral equations on planar
  curves with corners}.
\newblock \emph{Applied and Computational Harmonic Analysis}, 32:\penalty0
  45--64, 2012.

\bibitem[Camassa et~al.(2012)Camassa, McLaughlin, Moore, and Yu]{MooreJFM2012}
R.~Camassa, R.M. McLaughlin, M.N.J. Moore, and K.~Yu.
\newblock Stratified flows with vertical layering of density: experimental and
  theoretical study of flow configurations and their stability.
\newblock \emph{Journal of Fluid Mechanics}, 690:\penalty0 571--606, 2012.

\bibitem[Campbell et~al.(1996)Campbell, Ipsen, Kelley, Meyer, and
  Xue]{cam-ips-kel-mey-xue1996}
S.~L. Campbell, I.~C.~F. Ipsen, C.~T. Kelley, C.~D. Meyer, and Z.~Q. Xue.
\newblock {Convergence estimates for solution of integral equations with
  GMRES}.
\newblock \emph{Journal of Integral Equations and Applications}, 8\penalty0
  (1):\penalty0 19--34, 1996.

\bibitem[Cao et~al.(2010)Cao, Gunzburger, Hua, Wang, et~al.]{cao2010coupled}
Yanzhao Cao, Max Gunzburger, Fei Hua, Xiaoming Wang, et~al.
\newblock Coupled Stokes-Darcy model with Beavers-Joseph interface boundary
  condition.
\newblock \emph{Communications in Mathematical Sciences}, 8\penalty0
  (1):\penalty0 1--25, 2010.

\bibitem[Claudin et~al.(2017)Claudin, Dur{\'a}n, and
  Andreotti]{claudin2017dissolution}
Philippe Claudin, Orencio Dur{\'a}n, and Bruno Andreotti.
\newblock Dissolution instability and roughening transition.
\newblock \emph{Journal of Fluid Mechanics}, 832, 2017.

\bibitem[Cohen et~al.(2016)Cohen, Berhanu, Derr, and du~Pont]{cohen2016erosion}
Caroline Cohen, Michael Berhanu, Julien Derr, and Sylvain~Courrech du~Pont.
\newblock Erosion patterns on dissolving and melting bodies.
\newblock \emph{Physical Review Fluids}, 1\penalty0 (5):\penalty0 050508, 2016.

\bibitem[Cohen et~al.(2015)Cohen, Devauchelle, Seybold, Yi, Szymczak, and
  Rothman]{coh-dev-sey-yi-szy-rot2015}
Yossi Cohen, Olivier Devauchelle, Hansj{\"o}rg~F. Seybold, Robert~S. Yi, Piotr
  Szymczak, and Daniel~H. Rothman.
\newblock Path selection in the growth of rivers.
\newblock \emph{Proceedings of the National Academy of Sciences}, 112\penalty0
  (46):\penalty0 14132--14137, 2015.

\bibitem[Corona and Veerapaneni(2017)]{cor-vee2017}
Eduardo Corona and Shravan Veerapaneni.
\newblock {Boundary integral equation analysis for suspension of spheres in
  Stokes flow}.
\newblock \emph{arxiv}, 1707.06551, 2017.

\bibitem[Corona et~al.(2017)Corona, Greengard, Rachh, and
  Veerapaneni]{cor-gre-rac-vee2017}
Eduardo Corona, Leslie Greengard, Manas Rachh, and Shravan Veerapaneni.
\newblock {An integral equation formulation for rigid bodies in Stokes flow in
  three dimensions}.
\newblock \emph{Journal of Computational Physics}, 332:\penalty0 504--519,
  2017.

\bibitem[Coulier et~al.(2017)Coulier, Pouransari, and Darve]{cou-pou-dar2017}
Pieter Coulier, Hadi Pouransari, and Eric Darve.
\newblock {The Inverse Fast Multipole Method: Using a Fast Approximate Diret
  Solver as a Preconditioner for Dense Linear Systems}.
\newblock \emph{SIAM Journal on Scientific Computing}, 39\penalty0
  (3):\penalty0 A761--A796, 2017.

\bibitem[Daerr et~al.(2003)Daerr, Lee, Lanuza, and
  Cl{\'e}ment]{daerr2003erosion}
Adrian Daerr, Peter Lee, Jos{\'e} Lanuza, and {\'E}ric Cl{\'e}ment.
\newblock Erosion patterns in a sediment layer.
\newblock \emph{Physical Review E}, 67\penalty0 (6):\penalty0 065201, 2003.

\bibitem[Dallaston et~al.(2015)Dallaston, Hewitt, and
  Wells]{dallaston2015channelization}
MC~Dallaston, IJ~Hewitt, and AJ~Wells.
\newblock Channelization of plumes beneath ice shelves.
\newblock \emph{Journal of Fluid Mechanics}, 785:\penalty0 109--134, 2015.

\bibitem[de~Anna et~al.(2018)de~Anna, Quaife, Biros, and
  Juanes]{dea-qua-bir-jua2018}
Pietro de~Anna, Bryan Quaife, George Biros, and Ruben Juanes.
\newblock {Prediction of velocity distribution from pore structure in simple
  porous media}.
\newblock \emph{Physical Review Fluids}, 2\penalty0 (12):\penalty0 124103,
  2018.

\bibitem[Devauchelle et~al.(2012)Devauchelle, Petroff, Seybold, and
  Rothman]{Rothman2012}
O.~Devauchelle, A.P. Petroff, H.F. Seybold, and D.H. Rothman.
\newblock Ramification of stream networks.
\newblock \emph{Proceedings of the National Academy of Sciences}, 109\penalty0
  (51):\penalty0 20832--20836, 2012.

\bibitem[Domokos et~al.(2014)Domokos, Jerolmack, Sipos, and
  T{\"o}r{\"o}k]{domokos2014river}
Gabor Domokos, Douglas~J Jerolmack, Andras~{\'A} Sipos, and {\'A}kos
  T{\"o}r{\"o}k.
\newblock How river rocks round: resolving the shape-size paradox.
\newblock \emph{PloS One}, 9\penalty0 (2):\penalty0 e88657, 2014.

\bibitem[Gillman et~al.(2014)Gillman, Hao, and Martinsson]{gil-hao-mar2014}
A.~Gillman, S.~Hao, and P.G. Martinsson.
\newblock {A simplified technique for the efficient and highly accurate
  discretization of boundary integral equations in 2D on domains with corners}.
\newblock \emph{Journal of Computational Physics}, 256:\penalty0 214--219,
  2014.

\bibitem[Greenbaum et~al.(1992)Greenbaum, Greengard, and Mayo]{gre-gre-may1992}
Anne Greenbaum, Leslie Greengard, and Anita Mayo.
\newblock On the numerical solution of the biharmonic equation in the plane.
\newblock \emph{Physica D: Nonlinear Phenomena}, 60:\penalty0 216--225, 1992.

\bibitem[Greengard and Rokhlin(1987)]{gre-rok1987}
L.~Greengard and V.~Rokhlin.
\newblock {A Fast Algorithm for Particle Simulations}.
\newblock \emph{Journal of Computational Physics}, 73:\penalty0 325--348, 1987.

\bibitem[Greengard et~al.(1996)Greengard, Kropinski, and Mayo]{gre-kro-may1996}
L.~Greengard, M.C. Kropinski, and A.~Mayo.
\newblock {Integral Equation Methods for Stokes Flow and Isotropic Elasticity
  in the Plane}.
\newblock \emph{Journal of Computational Physics}, 125\penalty0 (2):\penalty0
  403--414, 1996.

\bibitem[Groen et~al.(2007)Groen, Gijsen, van~der Lugt, Ferguson, Hatsukami,
  van~der Steen, Yuan, and Wentzel]{gro-gij-van-fer-hat-van-yua-wen2007}
Harald~C. Groen, Frank~J.H. Gijsen, Aad van~der Lugt, Marina~S. Ferguson,
  Thomas~S. Hatsukami, Anton~F.W. van~der Steen, Chun Yuan, and Jolanda~J.
  Wentzel.
\newblock {Plaque Rupture in the Carotid Artery Is Localized at the High Shear
  Stress Region: A Case Report}.
\newblock \emph{Stroke}, 38:\penalty0 2379--2381, 2007.

\bibitem[Hanna(1969)]{han1969}
Steven~R. Hanna.
\newblock {The Formation of Longitudinal Sand Dunes by Large Helical Eddies in
  the Atmosphere}.
\newblock \emph{Journal of Applied Meteorology}, 8\penalty0 (6):\penalty0
  874--883, 1969.

\bibitem[Helsing(2011)]{hel2011}
J.~Helsing.
\newblock {A Fast and Stable Solver for Singular Integral Equations and
  Piecewise Smooth Curves}.
\newblock \emph{SIAM Journal on Scientific Comuting}, 33:\penalty0 153--174,
  2011.

\bibitem[Helsing and Ojala(2008)]{hel-oja2008a}
Johan Helsing and Rikard Ojala.
\newblock On the evaluation of layer potentials close to their sources.
\newblock \emph{Journal of Computational Physics}, 227:\penalty0 2899--2921,
  2008.

\bibitem[Hewett and Sellier(2017{\natexlab{a}})]{hewett2017evolution}
James~N Hewett and Mathieu Sellier.
\newblock Evolution of an eroding cylinder in single and lattice arrangements.
\newblock \emph{Journal of Fluids and Structures}, 70:\penalty0 295--313,
  2017{\natexlab{a}}.

\bibitem[Hewett and Sellier(2017{\natexlab{b}})]{hewett2017pear}
James~N Hewett and Mathieu Sellier.
\newblock The pear-shaped fate of an ice melting front.
\newblock \emph{arXiv preprint arXiv:1705.02536}, 2017{\natexlab{b}}.

\bibitem[Hou et~al.(1994)Hou, Lowengrub, and Shelley]{hou-low-she1994}
Thomas~Y. Hou, John~S. Lowengrub, and Michael~J. Shelley.
\newblock {Removing the Stiffness for Interfacial Flows with Surface Tension}.
\newblock \emph{Journal of Computational Physics}, 114:\penalty0 312--338,
  1994.

\bibitem[Huang et~al.(2015)Huang, Moore, and Ristroph]{Huang2015}
J.~M. Huang, M.N.J. Moore, and L.~Ristroph.
\newblock Shape dynamics and scaling laws for a body dissolving in fluid flow.
\newblock \emph{Journal of Fluid Mechanics}, 765, 2 2015.
\newblock ISSN 1469-7645.

\bibitem[Jerolmack et~al.(2012)Jerolmack, Ewing, Falcini, Martin, Masteller,
  Phillips, Reitz, and Buynevich]{jerolmack2012internal}
Douglas~J Jerolmack, Ryan~C Ewing, Federico Falcini, Raleigh~L Martin, Claire
  Masteller, Colin Phillips, Meredith~D Reitz, and Ilya Buynevich.
\newblock Internal boundary layer model for the evolution of desert dune
  fields.
\newblock \emph{Nature Geoscience}, 5\penalty0 (3):\penalty0 206, 2012.

\bibitem[Kabacao\u{g}lu et~al.(2018)Kabacao\u{g}lu, Quaife, and
  Biros]{kab-qua-bir2018}
Gokberk Kabacao\u{g}lu, Bryan Quaife, and George Biros.
\newblock {Low-resolution simulations of vesicle suspensions in 2D}.
\newblock \emph{Journal of Computational Physics}, 357:\penalty0 43--77, 2018.

\bibitem[Kl\"{o}ckner et~al.(2013)Kl\"{o}ckner, Barnett, Greengard, and
  O'Neil]{klo-bar-gre-one2013}
Andreas Kl\"{o}ckner, Alexander Barnett, Leslie Greengard, and Michael O'Neil.
\newblock {Quadrature by expansion: A new method for the evaluation of layer
  potentials}.
\newblock \emph{Journal of Computational Physics}, 252:\penalty0 332--349,
  2013.

\bibitem[Kondratiuk and Szymczak(2015)]{kondratiuk2015steadily}
Pawe{\l} Kondratiuk and Piotr Szymczak.
\newblock {Steadily Translating Parabolic Dissolution Fingers}.
\newblock \emph{SIAM Journal on Applied Mathematics}, 75\penalty0 (5):\penalty0
  2193--2213, 2015.

\bibitem[Ladyzhenskaya(1963)]{lad1963}
O.~A. Ladyzhenskaya.
\newblock \emph{{The Mathematical Theory of Viscous Incompressible Flow}}.
\newblock Gordon and Breach Science, 1963.

\bibitem[L{\'o}pez et~al.(2018)L{\'o}pez, Stickland, and
  Dempster]{lopez2018cfd}
Alejandro L{\'o}pez, Matthew~T Stickland, and William~M Dempster.
\newblock Cfd study of fluid flow changes with erosion.
\newblock \emph{Computer Physics Communications}, 2018.

\bibitem[Lu et~al.(2017)Lu, Rahimian, and Zorin]{lu-rah-zor2017}
Libin Lu, Abtin Rahimian, and Denis Zorin.
\newblock {Contact-aware simulations of particulate Stokesian suspensions}.
\newblock \emph{Journal of Computational Physics}, 347:\penalty0 160--182,
  2017.

\bibitem[Mitchell and Spagnolie(2016)]{mit-spa2016}
W.H. Mitchell and S.E. Spagnolie.
\newblock A generalized traction integral equation for {S}tokes flow, with
  applications to near-wall particle mobility and viscous erosion.
\newblock \emph{Journal of Computational Physics}, 2016.

\bibitem[Moore(2017)]{moore2017riemann}
M~Nicholas~J Moore.
\newblock {Riemann-Hilbert Problems for the Shapes Formed by Bodies Dissolving,
  Melting, and Eroding in Fluid Flows}.
\newblock \emph{Communications on Pure and Applied Mathematics}, 70\penalty0
  (9):\penalty0 1810--1831, 2017.

\bibitem[Moore et~al.(2013)Moore, Ristroph, Childress, Zhang, and
  Shelley]{moo-ris-chi-zha-she2013}
Matthew N.~J. Moore, Leif Ristroph, Stephen Childress, Jun Zhang, and
  Michael~J. Shelley.
\newblock Self-similar evolution of a body eroding in a fluid flow.
\newblock \emph{Physics of Fluids}, 25:\penalty0 116602, 2013.

\bibitem[Nienhuis et~al.(2014)Nienhuis, Perron, Kao, and
  Myrow]{nienhuis2014wavelength}
Jaap~H Nienhuis, J~Taylor Perron, Justin~CT Kao, and Paul~M Myrow.
\newblock Wavelength selection and symmetry breaking in orbital wave ripples.
\newblock \emph{Journal of Geophysical Research: Earth Surface}, 119\penalty0
  (10):\penalty0 2239--2257, 2014.

\bibitem[O.Pironneau(1973)]{pir1973}
O.Pironneau.
\newblock {On optimum profiles in Stokes flow}.
\newblock \emph{Journal of Fluid Mehanics}, 59\penalty0 (1):\penalty0 117--128,
  1973.

\bibitem[Perkins et~al.(2015)Perkins, Finnegan, and
  De~Silva]{perkins2015amplification}
Jonathan~P Perkins, Noah~J Finnegan, and Shanaka~L De~Silva.
\newblock Amplification of bedrock canyon incision by wind.
\newblock \emph{Nature Geoscience}, 8\penalty0 (4):\penalty0 305, 2015.

\bibitem[Perron et~al.(2009)Perron, Kirchner, and
  Dietrich]{perron2009formation}
J~Taylor Perron, James~W Kirchner, and William~E Dietrich.
\newblock Formation of evenly spaced ridges and valleys.
\newblock \emph{Nature}, 460\penalty0 (7254):\penalty0 502, 2009.

\bibitem[Petroff et~al.(2011)Petroff, Devauchelle, Abrams, Lobkovsky, Kudrolli,
  and Rothman]{petroff2011geometry}
Alexander~P Petroff, Olivier Devauchelle, Daniel~M Abrams, Alexander~E
  Lobkovsky, Arshad Kudrolli, and Daniel~H Rothman.
\newblock Geometry of valley growth.
\newblock \emph{Journal of Fluid Mechanics}, 673:\penalty0 245--254, 2011.

\bibitem[Picioreanu et~al.(2001)Picioreanu, van Loosdrecht, and
  Heijnen]{pic-van-hei2000}
Cristian Picioreanu, Mark C.~M. van Loosdrecht, and Joseph~J. Heijnen.
\newblock {Two-Dimensional Model of Biofilm Detachment Caused by Internal
  Stress from Liquid Flow}.
\newblock \emph{Biotechnology and Bioengineering}, 72\penalty0 (2):\penalty0
  205--218, 2001.

\bibitem[Power(1993)]{pow1993}
H.~Power.
\newblock {The completed double layer boundary integral equation method for
  two-dimensional Stokes flow}.
\newblock \emph{IMA Journal of Applied Mathematics}, 51\penalty0 (2):\penalty0
  123--145, 1993.

\bibitem[Power and Miranda(1987)]{pow-mir1987}
H.~Power and G.~Miranda.
\newblock Second kind integral equation formulation of Stokes' flows past a
  particle of arbitrary shape.
\newblock \emph{SIAM Journal on Applied Mathematics}, 47\penalty0 (4):\penalty0
  689--698, 1987.

\bibitem[Pozrikidis(1992)]{poz1992}
C.~Pozrikidis.
\newblock \emph{{Boundary Integral and Singularity Methods for Linearized
  Viscous Flow}}.
\newblock Cambridge University Press, New York, NY, USA, 1992.

\bibitem[Pozrikidis(1997)]{poz1997}
C.~Pozrikidis.
\newblock \emph{Introduction to Theoretical and Computational Fluid Dynamics}.
\newblock Oxford University Press, New York, 1997.

\bibitem[Quaife and Biros(2014)]{qua-bir2014a}
Bryan Quaife and George Biros.
\newblock High-volume fraction simulations of two-dimensional vesicle
  suspensions.
\newblock \emph{Journal of Computational Physics}, 274:\penalty0 245--267,
  2014.

\bibitem[Quaife and Biros(2015)]{qua-bir2015a}
Bryan Quaife and George Biros.
\newblock {On preconditioners for the Laplace double-layer potential in 2D}.
\newblock \emph{Numerical Linear Algebra with Applications}, 22:\penalty0
  101--122, 2015.

\bibitem[Quaife and Biros(2016)]{qua-bir2016}
Bryan Quaife and George Biros.
\newblock {Adaptive time stepping for vesicle simulations}.
\newblock \emph{Journal of Computational Physics}, 306:\penalty0 478--499,
  2016.

\bibitem[Quaife et~al.(2018)Quaife, Coulier, and Darve]{qua-cou-dar2018}
Bryan Quaife, Pieter Coulier, and Eric Darve.
\newblock {An efficient preconditioner for the fast simulation of a 2D Stokes
  flow in porous media}.
\newblock \emph{International Journal for Numerical Methods in Engineering},
  113:\penalty0 561--580, 2018.

\bibitem[Rachh and Serkh(2017)]{rac-ser2017}
M.~Rachh and K.~Serkh.
\newblock {On the solution of Stokes equation on regions with corners}.
\newblock \emph{arxiv}, 1711.04072, 2017.

\bibitem[Rachh and Greengard(2016)]{rac-gre2016}
Manas Rachh and L.~Greengard.
\newblock Integral equation methods for elastance and mobility problems in two
  dimensions.
\newblock \emph{SIAM Journal on Numerical Analysis}, 54\penalty0 (5):\penalty0
  2889--2909, 2016.

\bibitem[Ristroph(2018)]{ristroph2018sculpting}
Leif Ristroph.
\newblock Sculpting with flow.
\newblock \emph{Journal of Fluid Mechanics}, 838:\penalty0 1--4, 2018.

\bibitem[Ristroph et~al.(2012)Ristroph, Moore, Childress, Shelley, and
  Zhang]{ris-moo-chi-she-zha2012}
Leif Ristroph, Matthew N.~J. Moore, Stephen Childress, Michael~J. Shelley, and
  Jun Zhang.
\newblock Sculpting of an erodible body in flowing water.
\newblock \emph{Proceedings of the National Academy of Sciences}, 109\penalty0
  (48):\penalty0 19606--19609, 2012.

\bibitem[Rycroft and Bazant(2016)]{rycroft2016asymmetric}
Chris~H. Rycroft and Martin~Z. Bazant.
\newblock Asymmetric collapse by dissolution or melting in a uniform flow.
\newblock In \emph{Proc. R. Soc. A}, volume 472, page 20150531. The Royal
  Society, 2016.

\bibitem[Saad and Schultz(1986)]{saa-sch1986}
Y.~Saad and M.H. Schultz.
\newblock {GMRES: A Generalized Minimal Residual Algorithm for Solving
  Nonsymmetric Linear Systems}.
\newblock \emph{SIAM Journal on Scientific and Statistical Computing},
  7:\penalty0 856--869, 1986.

\bibitem[Sangani and Mo(1994)]{san-mo1994}
Ashok~S. Sangani and Guobiao Mo.
\newblock Inclusion of lubrication forces in dynamic simulations.
\newblock \emph{Physics of fluids}, 6\penalty0 (5):\penalty0 1653--1662, 1994.

\bibitem[Serkh and Rokhlin(2016)]{ser-rok2016}
K.~Serkh and V.~Rokhlin.
\newblock On the solution of elliptic partial differential equations on regions
  with corners.
\newblock \emph{Journal of Computational Physics}, 305:\penalty0 150--171,
  2016.

\bibitem[Shah(2002)]{sha2002}
Prediman~K. Shah.
\newblock {Pathophysiology of coronary thrombosis: role of plaque rupture and
  plaque erosion}.
\newblock \emph{Progress in Cardiovascular Diseases}, 44\penalty0 (5):\penalty0
  357--368, 2002.

\bibitem[Sidi and Israeli(1988)]{sid-isr1988}
Avram Sidi and Moshe Israeli.
\newblock Quadrature {M}ethods for {P}eriodic {S}ingular and {W}eakly
  {S}ingular {F}redholm {I}ntegral {E}quations.
\newblock \emph{Journal of Scientific Computing}, 3\penalty0 (2):\penalty0
  201--231, 1988.

\bibitem[Trefethen and Weideman(2014)]{tre-wei2014}
Lloyd~N. Trefethen and J.~A.~C. Weideman.
\newblock {The Exponentially Convergent Trapezoidal Rule}.
\newblock \emph{SIAM Review}, 56\penalty0 (3):\penalty0 385--458, 2014.

\bibitem[Ying et~al.(2004)Ying, Biros, and Zorin]{yin-bir-zor2004}
L.~Ying, G.~Biros, and D.~Zorin.
\newblock {A kernel-independent adaptive fast multipole algorithm in two and
  three dimensions}.
\newblock \emph{Journal of Computational Physics}, 194\penalty0 (2):\penalty0
  591--626, 2004.

\end{thebibliography}
